\documentclass[11pt]{amsart}
\usepackage{latexsym}
\usepackage{amssymb}
\usepackage{amsmath,amstext,amscd}
\usepackage{mathrsfs}
\usepackage{graphicx}
\usepackage{verbatim}
\usepackage[all]{xy}
\usepackage{amsfonts}
\usepackage{indentfirst}
\usepackage{enumerate}
\usepackage{txfonts}
\usepackage{pxfonts}
\usepackage{bezier,longtable}
\usepackage[labelformat=empty]{caption}

\newtheorem{theo}{Theorem}[section]

\newtheorem{lemma}[theo]{Lemma}

\newtheorem{corollary}[theo]{Corollary}
\newtheorem{prop}[theo]{Proposition}

\newtheorem{conjecture}[theo]{Conjecture}
\renewenvironment{proof}{ \emph{Proof}}{$\Box$}
\newtheorem{defi}{Definition}[section]

\newcommand {\ZZ} {\mathbb {Z}}

\newcommand {\CC} {\mathbb {C}}

\renewcommand{\ss}{\mathfrak{s}}

\renewcommand{\tt}{\mathfrak{t}}
\newcommand{\kk}{\mathfrak{k}}

\newcommand{\nn}{\mathfrak{n}}
\newcommand{\mm}{\mathfrak{m}}
\renewcommand{\aa}{\mathfrak{a}}
\newcommand{\hh}{\mathfrak{h}}

\newcommand{\qq}{\mathfrak{q}}
\newcommand{\pp}{\mathfrak{p}}
\renewcommand{\gg}{\mathfrak{g}}
\newcommand{\bb}{\mathfrak{b}}

\newcommand{\uu}{\mathfrak{u}}

\renewcommand{\ll}{\mathfrak{l}}

\renewcommand {\phi} {\varphi}

\newcommand{\half}{\frac{1}{2}}

\newcommand{\refth}[1]{Theorem \ref{#1}}
\newcommand{\refle}[1]{Lemma \ref{#1}}

\newcommand{\refcor}[1]{Corollary \ref{#1}}
\newcommand{\refprop}[1]{Proposition \ref{#1}}
\newcommand{\refcon}[1]{Conjecture \ref{#1}}
\newcommand{\refeq}[1]{(\ref{#1})}

\newcommand{\supp}{\mathrm{supp}}

\newcommand{\ad}{\mathrm{ad~}}

\newcommand{\ch}{\mathrm{ch}}
\newcommand{\Soc}{\mathrm{Soc}}
\newcommand{\pro}{\mathrm{pro}}
\newcommand{\ind}{\mathrm{ind}}
\newcommand{\Hom}{\mathrm{Hom}}

\renewcommand{\Top}{\mathrm{Top}}
\newcommand{\Ext}{\mathrm{Ext}}
\newcommand{\rk}{\mathrm{rk~}}
\newcommand{\field}[1]{\mathbb{#1}}

\def\cplus{\hbox{$\subset${\raise0.3ex\hbox{\kern -0.55em ${\scriptscriptstyle +}$}}\ }}

\def\clplus{\hbox{$\subset${\raise0.3ex\hbox{\kern -0.55em ${\scriptscriptstyle +}$}}\ }}
\def\crplus{\hbox{$\supset${\raise1.05pt\hbox{\kern -0.55em ${\scriptscriptstyle +}$}}\ }}

\title[On the structure of the fundamental series of generalized Harish-Chandra modules]{On the structure of the fundamental series of generalized Harish-Chandra modules}
\author[Ivan Penkov]{\;Ivan~Penkov}

\address{
Ivan Penkov
\newline School of Engineering and Science
\newline Jacobs University Bremen
\newline Campus Ring 1
\newline 28759 Bremen, Germany}
\email{i.penkov@jacobs-university.de}

\author[Gregg Zuckerman]{\;Gregg Zuckerman}

\address{
Gregg Zuckerman
\newline Department of Mathematics
\newline Yale University
\newline 10 Hillhouse Avenue, P.O. Box 208283
\newline New Haven, CT 06520-8283, USA}
\email{gregg@math.yale.edu}
\begin{document}
\maketitle
\tableofcontents

\newpage
\begin{abstract}
We continue the study of the fundamental series of generalized Harish-Chandra modules initiated in \cite{PZ2}. Generalized Harish-Chandra modules are $(\gg,\kk)$-modules of finite type where $\gg$ is a semisimple Lie algebra and $\kk\subset\gg$  is a reductive in $\gg$ subalgebra. A first result of the present paper is that a fundamental series module is a $\gg$-module of finite length. We then define the notions of strongly and weakly reconstructible simple $(\gg,\kk)$-modules $M$ which reflect to what extent $M$ can be determined via its appearance in the socle of a fundamental series module. 

In the second part of the paper we concentrate on the case $\kk\simeq sl(2)$ and prove a sufficient condition for strong reconstructibility. This strengthens our main result from \cite{PZ2} for the case $\kk =sl(2)$. We also compute the $sl(2)$-characters of all simple strongly reconstructible (and some weakly reconstructible) $(\gg,sl(2))$-modules. We conclude the paper by discussing a functor between a generalization of the category $\mathcal{O}$ and a category of $(\gg,sl(2))$-modules, and we conjecture that this functor is an equivalence of categories.

\end{abstract}
\section*{Introduction}
This paper is a continuation of our work \cite{PZ2}. By $\gg$ we denote a semisimple Lie algebra and by $\kk$ an arbitrary reductive in $\gg$ subalgebra. In \cite{PZ2} we introduced the fundamental series of generalized Harish-Chandra modules (or equivalently, $(\gg,\kk)$-modules of finite type over $\kk$) and proved that any simple generalized Harish-Chandra module with generic minimal $\kk$-type arises as the socle of an appropriate fundamental series module. Using this result we were able to show that any simple generalized Harish-Chandra module with generic minimal $\kk$-type can be reconstructed from its $\nn$-cohomology. This led to a classification of generalized Harish-Chandra modules with generic minimal $\kk$-type.

In the present paper we study the fundamental series further. After recalling the necessary preliminaries from \cite{PZ2}, we prove in Section 2 that any fundamental series generalized Harish-Chandra module has finite length. In Section 4 we introduce the concepts of a strongly reconstructible  and a weakly reconstructible simple generalized Harish-Chandra module. Theorem 3 from \cite{PZ2} implies that any simple Harish-Chandra module with generic minimal $\kk$-type is strongly reconstructible; however our aim is to study strong and weak reconstructibility of simple generalized Harish-Chandra modules which do not necessarily have a generic minimal $\kk$-type.

From Section 5 on, we concentrate on the case when the subalgebra $\kk$ of $\gg$ is isomorphic to $sl(2)$, i.e. we consider generalized Harish-Chandra $(\gg,sl(2))$-modules. Under this assumption we are able to considerably strengthen the results of \cite{PZ2} and establish new results about strong and weak reconstructibility. In particular, we prove that if $M$ is a simple $(\gg,sl(2))$-module with minimal $\kk$-type $V(\mu)$ satisfying $\mu\geq\frac{1}{2}(\lambda_1+\lambda_2)$ (note that $\mu\in\field{Z}_{\geq 0}$ as $\kk\simeq sl(2)$), $\lambda_1, \lambda_2$ being the maximum and submaximum eigenvalues in $\gg$ of a Cartan subalgebra $\tt$ of $\kk=sl(2)$, then $M$ is reconstructible by its $\nn$-cohomology. This yields a classification of simple $(\gg,sl(2))$-modules $M$ with $\mu\geq\frac{1}{2}(\lambda_1+\lambda_2)$ and proves that all such simple $(\gg,\kk)$-modules have finite type over $\kk$. For the principal $sl(2)$-subalgebra the bound $\frac{1}{2}(\lambda_1+\lambda_2)$ is linear in $\rk\gg$, while the bound established in \cite{PZ2} is cubic in $\rk\gg$. In addition, when $\kk$ is a direct summand of a symmetric subalgebra $\tilde{\kk}$ of $\gg$, we obtain new reconstruction results for Harish-Chandra modules.

In Section 7 we prove that for $\mu\geq\frac{1}{2}\lambda_1$, the socle of the fundamental series module is isomorphic to $R^1\Gamma_{\kk,\tt}(L_\pp(E))$, where $L_\pp(E)$ is the simple lowest weight module associated to the data $(\pp,E)$. The relative Kazhdan-Lusztig theory \cite{CC} yields an explicit formula for the $\tt$-character of $L_\pp(E)$. In turn, the theory of the derived Zuckerman functors yields an explicit formula for the $\kk$-character of the strongly reconstructible module  $R^1\Gamma_{\kk,\tt}(L_\pp(E))$. 

Section 8 is devoted to examples. We consider six explicit pairs $(\gg,sl(2))$ with $\rk\gg=2$ and we compute the respective sharp bounds on $\mu$ which ensure that a simple $(\gg,sl(2))$-module with minimal $sl(2)$-type $V(\mu)$ is strongly reconstructible. For a principal $sl(2)$-subalgebra of $sp(4)$ this sharp bound is the one established in the present paper, in the other five cases it turns out to be the bound from \cite{PZ2}.

In the final Section 9 we  discuss the possibility that, for $\kk\simeq sl(2)$ and a large enough $n$, the functor $R^1\Gamma_{\kk,\tt}$ is an equivalence between a certain category of $\bar{\pp}$-finite modules $\mathcal{C}_{\bar{\pp},\tt,n+2}$ and a category of $(\gg,\kk)$-modules $\mathcal{C}_{\kk,n}$. Proving or disproving this statement is an open problem. We conjecture that if $n\geq\frac{1}{2}(\lambda_1+\lambda_2)$, then $R^1\Gamma_{\kk,\tt}$ is an equivalence of categories between $\mathcal{C}_{\bar{\pp},\tt,n+2}$ and $\mathcal{C}_{\kk,n}$.

\textbf{Acknowledgements.} We thank Vera Serganova for pointing out to us Examples 1 and 2 in Subsection 8.7. Both authors acknowledge partial support through DFG Grant PE 980/3-1(SPP 1388). I. Penkov acknowledges the hospitality and partial support of Yale university, and G. Zuckerman acknowledges the hospitality of Jacobs University Bremen.

\section{Notation and preliminary results}
We start by recalling the setup of \cite{PZ2}.
\subsection{Conventions}

The ground field is $\field{C}$, and if not explicitly stated otherwise, all vector spaces and Lie algebras are defined over $\field{C}$. The sign $\otimes$ denotes tensor product over $\field{C}$. The superscript $^*$ indicates dual space. The sign $\clplus$ stands for semidirect sum of Lie algebras (if $\ll=\ll '\clplus\ll ''$, $\ll '$ is an ideal in $\ll$ and $\ll ''\cong\ll /\ll '$). $H^\cdot(\ll,M)$ stands for the cohomology of a Lie algebra $\ll$ with coefficients in an $\ll$-module $M$, and $M^\ll=H^0(\ll,M)$ stands for space of $\ll$-invariants of $M$. By $Z(\ll)$ we denote the center of $\ll$. $\Lambda^\cdot( \;)$ and $S^\cdot( \;)$ denote respectively the exterior and symmetric algebra.

If $\ll$ is a Lie algebra, then $U(\ll)$ stands for the enveloping algebra of $\ll$ and $Z_{U(\ll)}$ denotes the center of $U(\ll)$. We identify $\ll$-modules with $U(\ll)$-modules. It is well known that if $\ll$ is finite dimensional and $M$ is a simple $\ll$-module (or equivalently a simple $U(\ll)$-module), $Z_{U(\ll)}$ acts on $M$ via a $Z_{U(\ll)}$-\textit{character}, i.e. via an algebra homomorphism $\theta_{M}: Z_{U(\ll)}\rightarrow\field{C}$.

We say that an $\ll$-module $M$ is \textit{generated} by a subspace $M'\subset M$ if $U(\ll)\cdot M'=M$, and we say that $M$ is \textit{cogenerated} by $M'\subset M$, if for any non-zero homomorphism $\psi :M\rightarrow\bar{M}$, $M'\cap\ker\psi=\{0\}$. By $\Soc M$ we denote the socle (i.e. the unique maximal semisimple submodule) of an $\ll$-module $M$; by $\Top M$ we denote the unique maximal semisimple quotient of $M$, when $M$ has finite length.

If $\ll$ is a Lie algebra, $M$ is an $\ll$-module, and $\omega\in\ll^*$, we put $M^\omega:=\{m\in M\; |\; \ll\cdot m=\omega(\ll)m \;\forall l \in\ll\}$. By $\supp_\ll M$ we denote the set $\{\omega\in\ll^*\; |\; M^\omega\neq 0\}$.

A finite \textit{multiset} is a function $f$ from a finite set $D$ into $\field{N}$. A \textit{submultiset} of $f$ is a multiset $f'$ defined on the same domain $D$ such that $f'(d)\leq f(d)$ for any $d\in D$. For any finite multiset $f$, defined on an additive monoid $D$, we can put $\rho_f:=\frac{1}{2}\sum_{d\in D}f(d)d$.

If $\dim M<\infty$ and $M=\bigoplus_{\omega\in\ll^*}M^\omega$, then $M$ determines the finite multiset $\ch_{\ll}M$ which is the function $\omega\mapsto\dim M^\omega$ defined on $\supp_{\ll}M$.

\subsection{Reductive subalgebras, compatible parabolics and generic $\kk$-types}

Let $\gg$ be a finite-dimensional semisimple Lie algebra. By $\gg$-mod we denote the category of $\gg$-modules. Let $\kk\subset\gg$ be an algebraic subalgebra which is reductive in $\gg$. We fix a Cartan subalgebra $\tt$ of $\kk$ and a Cartan subalgebra $\hh$ of $\gg$ such that $\tt\subset\hh$. By $\Delta$ we denote the set of $\hh$-roots of $\gg$, i.e. $\Delta=\{\supp_\hh\gg\}\setminus\{0\}$. Note that, since $\kk$ is reductive in $\gg$, $\gg$ is a $\tt$-weight module, i.e. $\gg=\bigoplus_{\lambda\in\tt^*}\gg^\lambda$. We set $\Delta_\tt:=\{\supp_\tt\gg\}\setminus\{0\}$. Note also that the $\field{R}$-span of the roots of $\hh$ in $\gg$ fixes a real structure on $\hh^*$, whose projection onto $\tt^*$ is a well-defined real structure on $\tt^*$. In what follows, we will denote by $\mathrm{Re}\lambda$ the real part of an element $\lambda\in\tt^*$. We fix also a Borel subalgebra $\bb_\kk\subset\kk$ with $\bb_\kk\supset\tt$. Then $\bb_\kk=\tt\crplus\nn_\kk$, where $\nn_\kk$ is the nilradical of $\bb_\kk$. We set $\rho:=\rho_{\ch_\tt\nn_\kk}$. The quintet $\gg,\hh,\kk,\bb_\kk,\tt$ will be fixed throughout the paper. By $W$ we denote the Weyl group of $\gg$, and by $C(\cdot)$ - centralizer in $\gg$.

As usual, we will parametrize the characters of $Z_{U(\gg)}$ via the Harish-Chandra homomorphism. More precisely, if $\bb$ is a given Borel subalgebra of $\gg$ with $\bb\supset\hh$ ($\bb$ will be specified below), the $Z_{U(\gg)}$-character corresponding to $\kappa\in\hh^*$ via the Harish-Chandra homomorphism defined by $\bb$ will be denoted by $\theta_\kappa$ ($\theta_{\rho_{\ch_\hh\bb}}$ is the trivial $Z_{U(\gg)}$-character).

By $\langle\;,\;\rangle$ we denote the unique $\gg$-invariant symmetric bilinear form on $\gg^*$ such that $\langle\alpha,\alpha\rangle=2$ for any long root of a simple component of $\gg$. The form $\langle\;,\;\rangle$ enables us to identify $\gg$ with $\gg^*$. Then $\hh$ is identified with $\hh^*$, and $\kk$ is identified with $\kk^*$. We will sometimes consider $\langle\;,\;\rangle$ as a form on $\gg$. The superscript $\perp$ indicates orthogonal space. Note that there is a canonical $\kk$-module decomposition $\gg=\kk\oplus\kk^\perp$. We also set $\parallel\kappa\parallel^2:=\langle\kappa,\kappa\rangle$ for any $\kappa\in\hh^*$.

We say that an element $\lambda\in\tt^*$ is $(\gg,\kk)$-\textit{regular} if $\langle\mathrm{Re}\lambda,\sigma\rangle\neq 0$ for all $\sigma\in\Delta_\tt$. To any $\lambda\in\tt^*$ we associate the following parabolic subalgebra $\pp_\lambda$ of $\gg$: $$\pp_\lambda=\hh\oplus(\bigoplus_{\alpha\in\Delta_\lambda}\gg^\alpha),$$ where $\Delta_\lambda:=\{\alpha\in\Delta\; |\; \langle\mathrm{Re}\lambda,\sigma\rangle\geq 0\}$. By $\mm_\lambda$ and $\nn_\lambda$ we denote respectively the reductive part of $\pp$ (containing $\hh$) and the nilradical of $\pp$. In particular $\pp_\lambda=\mm_\lambda\crplus\nn_\lambda$, and if $\tt$ is $\bb_\kk$-dominant, then $\pp_\lambda\cap\kk=\bb_\lambda$. We call $\pp_\lambda$ a \textit{$\tt$-compatible parabolic subalgebra}. A $\tt$-compatible parabolic subalgebra $\pp=\mm\crplus\nn$ (i.e. $\pp=\pp_\lambda$ for some $\lambda\in\tt^*$) is \textit{minimal} if it does not properly contain another $\tt$-compatible parabolic subalgebra. It is an important observation that if $\pp=\mm\crplus\nn$ is minimal, then $\tt\subset Z(\mm)$. In fact, a $\tt$-compatible parabolic subalgebra $\pp$ is minimal if and only if $\mm$ equals the centralizer $C(\tt)$ of $\tt$ in $\gg$, or equivalently if and only if $\pp=\pp_\lambda$ with $\lambda$ $(\gg,\kk)$-regular. In this case $\nn\cap\kk =\nn_\kk$.

Any $\tt$-compatible parabolic subalgebra $\pp =\pp_\lambda$ has a well-defined opposite parabolic subalgebra $\bar{\pp}:=\pp_{-\lambda}$; clearly $\pp$ is minimal if and only if $\bar{\pp}$ is minimal.
\begin{lemma}
If $\mathrm{Re}\lambda([\kk,\kk]\cap\tt)\neq 0$, then $\pp_\lambda$ and $\kk$ generate the Lie algebra $\gg$.
\end{lemma}
\begin{proof}
Since Re$\lambda([\kk,\kk]\cap\tt)\neq 0$, there exists an $sl(2)$-subalgebra $\kk'\subset\kk$ such that Re$\lambda(\tt\cap\kk ')\neq 0$. By definition, $\pp_\lambda$ contains all $\kk'$-singular vectors of $\gg$. Hence $\pp_\lambda$ generates $\gg$ as a $\kk'$-module, i.e. $\pp_\lambda$ and $\kk'$ generate $\gg$.
\end{proof}

A \textit{$\kk$-type} is by definition a simple finite-dimensional $\kk$-module. By $V(\mu)$ we denote a $\kk$-type with $\bb_\kk$-highest weight $\mu$ ($\mu$ is then $\kk$-integral and $\bb_\kk$-dominant). Let $V(\mu)$ be a $\kk$-type such that $\mu+2\rho$ is $(\gg,\kk)$-regular, and let $\pp=\mm\crplus\nn$ be the minimal compatible parabolic subalgebra $\pp_{\mu+2\rho}$. Put $\tilde{\rho}_\nn :=\rho_{\ch_\hh \nn}$ and $\rho_\nn:=\rho_{\ch_\tt\nn}$. Clearly $\rho_\nn =\tilde{\rho}_\nn |_\tt$. We define $V(\mu)$ to be \textit{generic} if the following two conditions hold:
\begin{enumerate}
\item $\langle\mathrm{Re}\mu+2\rho-\rho_\nn,\alpha\rangle\geq 0\;\forall\alpha\in\supp_\tt\nn_\kk$;
\item $\langle\mathrm{Re}\mu+2\rho-\rho_S,\rho_S\rangle >0$ for every submultiset $S$ of $\ch_\tt\nn$.
\end{enumerate}

It is easy to show that there exists a positive constant $C$ depending only on $\gg,\kk$ and $\pp$ such that $\langle\mathrm{Re}\mu+2\rho,\alpha\rangle >C$ for every $\alpha\in\supp_\tt\nn$ implies $\pp_{\mu+2\rho}=\pp$ and that $V(\mu)$ is generic.

In agreement with \cite{PZ2}, we define a $\gg$-module $M$ to be a $(\gg,\kk)$-\textit{module} if $M$ is isomorphic as a $\kk$-module to a direct sum of isotypic components of $\kk$-types. If $M$ is a $(\gg,\kk)$-module, we write $M[\mu]$ for the $V(\mu)$-isotypic component of $M$, and we say that $V(\mu)$ is a $\kk$-\textit{type of} $M$ if $M[\mu]\neq 0$. We say that a $(\gg,\kk)$-module $M$ is \textit{of finite type} if $\dim M[\mu]\neq\infty$ for every $\kk$-type $V(\mu)$. We will also refer to $(\gg,\kk)$-modules of finite type as \textit{generalized Harish-Chandra modules}.

Let $\Theta_\kk$ be the discrete subgroup of $Z(\kk)^*$ generated by $\supp_{Z(\kk)}\gg$. By $\mathcal{M}$ we denote the class of $(\gg,\kk)$-modules $M$ for which there exists a finite subset $S\subset Z(\kk)^*$ such that $\supp_{Z(\kk)}M\subset(S+\Theta_\kk)$. If $M$ is a module in $\mathcal{M}$, a $\kk$-type $V(\mu)$ of $M$ is \textit{minimal} if the function $\mu '\mapsto \parallel\mathrm{Re}\mu '+2\rho \parallel^2 $ defined on the set $\{\mu '\in\tt^*\; |\; M[\mu ']\neq 0\}$ has a minimum at $\mu$. Any non-zero $(\gg,\kk)$-module $M$ in $\mathcal{M}$ has a minimal $\kk$-type. This follows from the fact that the squared length of a vector has a minimum on every shifted lattice in Euclidean space.

\subsection{The fundamental series of generalized Harish-Chandra modules}

Recall that the \textit{functor of $\kk$-locally finite vectors} $\Gamma_{\kk,\tt}$ is a well-defined left exact functor on the category of $(\gg,\tt)$-modules with values in $(\gg,\kk)$-modules, $$\Gamma_{\kk,\tt}(M)=\sum_{M'\subset M, \dim M'=1, \dim U(\kk)\cdot M'<\infty}M'.$$ By $R^\cdot\Gamma_{\kk,\tt}:=\bigoplus_{i\geq 0}R^i\Gamma_{\kk,\tt}$ we denote as usual the total right derived functor of $\Gamma_{\kk,\tt}$, see \cite{PZ1} and the references therein.

If $M$ is a $(\gg,\kk)$-module of finite type, then $\Gamma_{\kk,0}(M^*)$ is a well-defined $(\gg,\kk)$-module of finite type and $\Gamma_{\kk,0}(\cdot^*)$ is an involution on the category of $(\gg,\kk)$-modules of finite type. We put $\Gamma_{\kk,0}(M^*)=M_\kk^*$. There is an obvious $\gg$-invariant non-degenerate pairing $M\times M_\kk^*\rightarrow\field{C}$.

\begin{lemma}\label{lemma1.2}
Let $W$ be a finite-dimensional $\gg$-module and $M$ be a finite length $(\gg,\kk)$-module of finite type over $\kk$. Then

a) $W\otimes M$ is a $(\gg,\kk)$-module of finite type.

b) $W\otimes M$ is a $\gg$-module of finite length.
\end{lemma}
\begin{proof}
a) Since $\kk$ is finite dimensional and reductive in $\gg$, the class of $(\gg,\kk)$-modules is closed under tensor products. Let $V(\mu)$ be a $\kk$-type. Since $W$ is finite dimensional, $\Hom_\kk(V(\mu), W\otimes M)\cong\Hom_\kk(V(\mu)\otimes W^*,M)$, which is finite dimensional, since $V(\mu)\otimes W^*$ is finite dimensional and $M$ has finite type over $\kk$.

b) Since $M$ has finite length, $M$ is finitely generated over $\gg$. Note that $M_\kk^*$, the $\kk$-finite dual of $M$ is a $(\gg,\kk)$-module of finite length and hence $M_\kk^*$ is finitely generated over $\gg$ and likewise $W^*\otimes M_\kk^*$ is finitely generated. Hence, $(W\otimes M)_\kk^*$ is finitely generated, and satisfies the ascending chain condition. We have already seen that $W\otimes M$ is finitely generated; thus $W\otimes M$ satisfies the ascending chain condition. We conclude that $W\otimes M$ has finite length.
\end{proof}

We also introduce the following notation: if $\qq$ is a subalgebra of $\gg$ and $J$ is a $\qq$-module, we set $\ind^{\gg}_\qq J:=U(\gg)\otimes_{U(\qq)}J$ and $\pro_\qq^\gg J:=\Hom_{U(\qq)}(U(\gg),J)$. For a finite-dimensional $\pp$- or $\bar{\pp}$-module $E$ we set $N_\pp(E):=\Gamma_{\tt,0}(\pro_\pp^\gg(E\otimes\Lambda^{\dim \nn}(\nn)))$, $N_{\bar{\pp}}(E^*):=\Gamma_{\tt,0}(\pro_{\bar{\pp}}^\gg(E^*\otimes\Lambda^{\dim \nn}(\nn^*)))$. Note that both $N_\pp(E)$ and $N_{\bar{\pp}}(E^*)$ have simple socles, as long as $E$ itself is simple.

The \textit{fundamental series} of $(\gg,\kk)$-modules of finite type $F^\cdot(\pp,E)$ is defined as follows. Let $\pp=\mm\crplus\nn$ be a minimal compatible parabolic subalgebra, $E$ be a simple finite dimensional $\pp$-module on which $\tt$ acts via the weight $\omega\in\tt^*$, and $\mu:=\omega+2\rho_\nn^\perp$ where $\rho_\nn^\perp :=\rho_\nn -\rho$. Set $$F^\cdot(\pp,E):=R^\cdot\Gamma_{\kk,\tt}(N_\pp(E)).$$ 

Then the following assertions hold under the assumptions that $\pp=\pp_{\mu+2\rho}$ and that $\mu$ is $\bb_\kk$-dominant, $\kk$-integral and yields a generic $\kk$-type $V(\mu)$ (Theorem 2 of \cite{PZ2}).
\begin{enumerate}[a)]
\item $F^\cdot(\pp,E)$ is a $(\gg,\kk)$-module of finite type in the class $\mathcal{M}$, and $Z_{U(\gg)}$ acts on $F^\cdot(\pp,E)$ via the $Z_{U(\gg)}$-character $\theta_{\nu+\tilde{\rho}}$ where $\tilde{\rho}:=\rho_{\ch_\hh\bb}$ for some fixed Borel subalgebra $\bb$ of $\gg$ with $\bb\supset\hh,\;\bb\subset\pp$ and $\bb\cap\kk=\bb_\kk$, and where $\nu$ is the $\bb$-highest weight of $E$ (note that $\nu|_\tt=\omega$).
\item $F^i(\pp,E)=0$ for $i\neq s:=\dim\nn_\kk$ .
\item There is a $\kk$-module isomorphism $$F^s(\pp,E)[\mu]\cong\field{C}^{\dim E}\otimes V(\mu),$$ and $V(\mu)$ is the unique minimal $\kk$-type of $F^s(\pp,E)$.
\item Let $\bar{F}^s(\pp,E)$ be the $\gg$-submodule of $F^s(\pp,E)$ generated by $F^s(\pp,E)[\mu]$. Then $\bar{F}^s(\pp,E)$ is the unique simple submodule of $F^s(\pp,E)$, and moreover, $F^s(\pp,E)$ is cogenerated by $F^s(\pp,E)[\mu]$. This implies that $F^s(\pp,E)^*_\tt$ is generated by $F^s(\pp,E)^*_\tt[w_\mm(-\mu)]$, where $w_\mm\in W_\kk$  is the element of maximal length in the Weyl group $W_\kk$ of $\kk$.
\item For any non-zero $\gg$-submodule $M$ of $F^s(\pp,E)$ there is an isomorphism of $\mm$-modules $$H^r(\nn,M)^\omega\cong E.$$
\end{enumerate}

\section{On the fundamental series of $(\gg,\kk)$-modules}

In the rest of the paper, $\pp$ is a minimal $\tt$-compatible parabolic subalgebra and $E$ is a simple finite-dimensional $\pp$-module. Then $\nn\cdot E=0$ and $E$ is a simple $\mm=C(\tt)$-module. Fix a Borel subalgebra $\bb_{\mm}$ in $\mm$ such that $\hh\subseteq\bb_{\mm}$. Write $\bb=\bb_{\mm}\crplus\nn$; then $\bb$ is a Borel subalgebra of $\gg$. Set $\tilde{\rho} = \rho_{\ch_{\hh} \bb} $.
\begin{theo}\label{th2.1}
Assume that $\pp=\pp_{\mu+2\rho}$ and that $\mu$ is generic. Assume in addition that $N_\pp(E)$ is a simple $\gg$-module. Then $F^s(\pp,E)$ is a simple (in particular, non-zero) $\gg$-module.  
\end{theo}
\begin{proof}
 By the Duality Theorem from \cite{EW},
 \begin{equation}\label{eqn1}
 (R^i\Gamma_{\kk,\tt}(X))^*_\kk\cong R^{2s-i}\Gamma_{\kk,\tt}(X_\tt^*)
 \end{equation}
 for any $(\gg,\tt)$-module $X$ of finite type over $\tt$. Set $X=N_\pp(E)$. Then $X$ is a $(\gg,\tt)$-module of finite type over $\tt$ (see for instance \cite{Z}), and \refeq{eqn1} yields for $i=s$ $$F^s(\pp,E)_\kk^*\cong R^s\Gamma_{\kk,\tt}(N_\pp(E)_\tt^*).$$

We have $(\ind_\pp^\gg(E^*\otimes\Lambda^{\dim \nn}(\nn)^*))^*\cong\pro_\pp^\gg(E\otimes\Lambda^{\dim\nn}(\nn)).$ Thus $$N_\pp(E)=\Gamma_{\kk,0}(\pro_\pp^\gg(E\otimes\Lambda^{\dim\nn}(\nn)))\cong (\ind_\pp^\gg(E^*\otimes\Lambda^{\dim\nn}(\nn^*)))^*_\tt.$$ Moreover, $N_\pp(E)$ has finite type over $\tt$. Hence, $N_\pp(E)^*_\tt\cong\ind_\pp^\gg(E^*\otimes\Lambda^{\dim\nn}(\nn^*))$.
 
By Frobenius reciprocity, there is a canonical $\gg$-module homomorphism $$N_\pp(E)_\tt^*\cong\ind_\pp^\gg(E^*\otimes\Lambda^{\dim\nn}(\nn^*))\stackrel{\varphi}{\rightarrow} N_{\bar{\pp}}(E^*)$$ whose restriction to $E^*\otimes\Lambda^{\dim\nn}(\nn^*)$ is the identity. As $N_\pp(E)_\tt^*$ is simple by our assumption, $\varphi$ is injective. Moreover, $\varphi$ must be surjective as the $\tt$-characters of $N_\pp(E)_\tt^*$ and $N_{\bar{\pp}}(E^*)$ are equal. Therefore there is a commutative diagram of isomorphisms
 $$\begin{array}[c]{cccc}
& R^s\Gamma_{\kk,\tt}(N_\pp(E)^*_\tt)& \stackrel{\sim}{\rightarrow}& R^s\Gamma_{\kk,\tt}(N_{\bar{\pp}}(E^*))\\
 &\wr\downarrow&&\wr\downarrow\\
 \gamma :& F^s(\pp,E)_\kk^*& \stackrel{\sim}{\rightarrow}&F^s(\bar{\pp},E^*).
 \end{array}$$
 
 The fact that $V(\mu)$ is generic for $\pp$ implies immediately that $V(\mu)^*$ is generic for $\bar{\pp}$. Thus $F^s(\pp,E)$ is cogenerated by its minimal $\kk$-isotypic component $F^s(\pp,E)[\mu]$, and $F^s(\bar{\pp},E^*)$ cogenerated by its minimal $\kk$-isotypic component. On the other hand, the isomorphism $\gamma$ implies that $F^s(\bar{\pp},E^*)$ is also generated by its minimal $\kk$-isotypic component as $F^s(\pp,E)^*_\kk$ is generated by its minimal $\kk$-isotypic component. We conclude that $F^s(\bar{\pp},E^*)\cong F^s(\pp,E)^*_\kk$ is simple, which in turn shows that $F^s(\pp,E)\cong(F^s(\pp,E)^*_\kk)^*_\kk$ is simple.
\end{proof}

Assume that the $\bb$-highest weight of $E$ is $\nu\in\hh^*$. Set $\omega:=\nu|_\tt$ and $\mu:=\omega+2\rho_{\nn}^{\perp}$.

\begin{corollary}\label{cor2.3}
Let $\nu+\tilde{\rho}$ be $\bb$-dominant, i.e. $\mathrm{Re}\langle\nu+\tilde{\rho},\gamma\rangle\geq 0$ for any root $\gamma$ of $\hh$ in $\bb$. Then, under the assumption that $\mu$ is generic and that $\pp=\pp_{\mu+2\rho}$, we have $F^i(\pp,E)=0$ for $i\neq s$. Moreover $F^s(\pp,E)$ is simple. Thus, $F^i(\pp,E)$ has finite length for all $i\geq 0$.
\end{corollary}
\begin{proof}
Under the hypothesis on $\nu$, $\ind_\pp^\gg(E^*\otimes\Lambda^{\dim\nn}(\nn^*))\cong N_\pp(E)_\tt^*$ is simple, hence $N_\pp(E)$ is simple, and the statement follows from \refth{th2.1}
\end{proof}

In the rest of this section, $\pp$ is an arbitrary minimal $\tt$-compatible parabolic subalgebra.
\begin{lemma}\label{le2.4}
For any $C\in\field{Z}_{\geq 0}$, there exists a $\bb$-dominant integral weight $\sigma_{0}$ of $\gg$ such that $\langle\sigma_{0}|_\tt,\ \alpha\rangle > C$ for every weight $\alpha$ of $\tt$ in $\nn$.
\end{lemma}
\begin{proof}
Since $\pp$ is $\tt$-compatible, there exists $\kappa\in\tt^*$ such that $\frac{2\langle\kappa,\ \gamma|_\tt\rangle}{\langle\gamma,\ \gamma\rangle}\in\field{Z}_{\geq 0}$ for every root $\gamma$ of $\hh$ in $\nn$. Regard $\tt^*$ as a subspace of $\hh^*$ via the Killing form of $\gg$ restricted to $\hh$. Then, $\frac{2\langle\kappa,\ \gamma\rangle}{\langle\gamma,\ \gamma\rangle}\in\field{Z}_{\geq 0}$ for every root $\gamma$ of $\hh$ in $\nn$ and $\langle\kappa,\ \gamma\rangle=0$ for every root $\gamma$ of $\hh$ in $\mm$. Hence, $\kappa$ is a dominant integral weight of $\hh$ in $\gg$. Finally, choose a positive integer $r$ such that $\langle r\kappa,\ \alpha\rangle > C$ for every weight $\alpha$ of $\tt$ in $\nn$. Then $\sigma_{0}=r\kappa$ is a $\bb$-dominant weight of $\gg$ as required.
\end{proof}

\begin{prop}\label{pr2.5}
Suppose that $\nu+\tilde{\rho}$ is $\bb$-dominant. Then, $F^i(\pp,E)=0$ for $i\neq s$, and $F^s(\pp,E)$ has finite length.
\end{prop}
\begin{proof}
Fix a constant $C\in\field{Z}_{>0}$. Chose $m\in\field{Z}_{>0}$ such that $\langle\mathrm{Re}\mu+2\rho,\alpha\rangle\geq-mC$ for every $\alpha\in\supp_\tt\nn$. Set $\sigma_1=(m+1)\sigma_0$, where $\sigma_0$ is defined in \refle{le2.4}. Then $\langle\mathrm{Re}\mu+\sigma_1+2\rho,\alpha\rangle\geq C$ for every $\alpha\in\supp_\tt\nn$, and by possibly making $C$ larger, we can assume that $\mu+\sigma_1$ is generic. In particular, $\pp_{\mu+\sigma_1+2\rho}=\pp$.

Next, let $\nu_1:=\nu+\sigma_1$, and $E_1$ be a simple finite-dimensional $\mm$-module with highest weight $\nu_1$. Set $\mu_1:=\mu+\sigma_1 |_\tt$. Then, by \refcor{cor2.3}, $F^i(\pp,E_1)=0$ for $i\neq s$, and $F^s(\pp,E_1)$ is a simple $\gg$-module. Furthermore, by Propositions 2.6 and 2.12 in \cite{Z}, $F^i(\pp,E)$ is a direct summand of $V_\gg(\sigma_1)^*\otimes F^i(\pp,E_1)$ where $V_\gg(\sigma_1)$ stands for the finite-dimensional $\gg$-module with $\bb$-highest weight $\sigma_1$. \refle{lemma1.2} implies the statement.
\end{proof}

\textbf{Remark.} By a more refined argument with translation functors \cite{BG} one can show using the result of \cite{PZ3} that $F^s(\pp,E)$ is simple and hence non-zero,while $F^{i}(\pp,E)=0$ if $i\neq s$.
\begin{theo}\label{th2.6}
The $(\gg,\kk)$-module $F^i(\pp,E)$ has finite length for any simple $\pp$-module $E$ and any $i\in\field{Z}_{\geq 0}$.
\end{theo}
\begin{proof}
We will assume at first that $\nu+\tilde{\rho}$ is a regular weight of $\hh$ in $\gg$. Then, there exists a unique element $w\in W$ such that $w^{-1}(\nu+\tilde{\rho})$ is dominant for $\hh$ in $\gg$. Denote by $d(\nu)$ the length $l(w)$. We will argue by induction on $d(\nu)$. The theorem is true for $d(\nu)=0$ by \refprop{pr2.5}.

Suppose we assume the theorem for $d(\nu)=n\in\field{Z}_{>0}$, $n$ being fixed. If $d(\nu)=n+1$, we can choose a root $\gamma$ of $\hh$ in $\gg$ such that $d(s_{\gamma}(\nu))=n$. Let $D$ have highest weight $s_{\gamma}(\nu+\tilde{\rho})-\tilde{\rho}$. We will show that the finiteness of the length of $F^i(\pp,D)$ for all $i$ implies the finiteness of the length of $F^i(\pp,E)$ for all $i$.

Case I: $\frac{2\langle\mathrm{Re}\nu+\tilde{\rho},\ \gamma\rangle}{\langle\gamma,\ \gamma\rangle}\in\field{Z}_{\geq 0}$.

Choose a translation functor $\Psi$ so that $\Psi(N_\pp(D))$ has a central character which is singular with respect to precisely $\gamma$. Let $\Phi$ be the translation functor adjoint to $\Psi$.
By highest weight module theory we have a short exact sequence $$0\rightarrow N_\pp(D)\rightarrow\Phi\circ\Psi(N_\pp(D))\rightarrow N_\pp(E)\rightarrow 0.$$
This short exact sequence yields a long exact sequence $$...\rightarrow R^i\Gamma_{\kk,\tt}(N_\pp(D))\rightarrow R^i\Gamma_{\kk,\tt}\Phi\circ\Psi(N_\pp(D))\rightarrow R^i\Gamma_{\kk,\tt}(N_\pp(E))\rightarrow R^{i+1}\Gamma_{\kk,\tt}(N_\pp(D))\rightarrow ...\; .$$
We can rewrite this long exact sequence as $$...\rightarrow\Phi\circ\Psi(F^i(\pp,D))\rightarrow F^i(\pp,E)\rightarrow F^{i+1}(\pp,D)\rightarrow ... .$$

By assumption, $F^i(\pp,D)$ and $F^{i+1}(\pp,D)$ have finite length. Hence, by \refle{lemma1.2} $\Phi\circ\Psi(F^i(\pp,D))$ has finite length. By the long exact sequence, $F^{i}(\pp,E)$ has finite length.

Case II: $\frac{2\langle\mathrm{Re}\nu+\tilde{\rho},\ \gamma\rangle}{\langle\gamma,\ \gamma\rangle}\notin\field{Z}_{\geq 0}$.

Choose an integral weight $\sigma\in\hh^*$ such that if $w(\nu+\tilde{\rho})$ is dominant, $s_{\gamma}w(\nu-\sigma+\tilde{\rho})$ is dominant. Let $D$ be a finite-dimensional simple $\mm$-module such that the highest weight of $D$ is $\nu-\sigma$. There exists a translation functor $\Psi$ such that $\Psi(N_\pp(E))=N_\pp(D)$. 

For the adjoint functor $\Phi$, $\Phi(N_\pp(D))=N_\pp(E)$. Then $R^{i}\Gamma_{\kk,\tt}(\Phi(N_\pp(D)))=R^{i}\Gamma_{\kk,\tt}(N_\pp(E))$. 

Hence, $F^{i}(\pp,E)=\Phi (F^{i}(\pp,D))$, and by \refle{lemma1.2} the finiteness of the length of $F^{i}(\pp,D)$ implies the same for $F^{i}(\pp,E)$. 
\end{proof}
\begin{corollary}
Let $A$ be a finite-dimensional $(\pp,\tt)$-module. Then $R^{i}\Gamma_{\kk,\tt}(N_\pp(A))$ has finite length for all $i$.
\end{corollary}
\begin{proof}
Induction on the length of $A$ as a $\pp$-module: if $A$ has length 1, then $\nn\cdot A=0$, and we are back to \refth{th2.6}.
\end{proof}

In what follows we denote by $\mathcal{C}_{\bar{\pp},\tt}$ the full subcategory of $\gg$-modules consisting of finitely generated $(\gg,\tt)$-modules which are locally $\bar{\pp}$-finite.
\begin{lemma}\label{le2.8}
Let $N\in\mathcal{C}_{\bar{\pp},\tt}$. Then $N$ has finite length, each simple constituent of $N$ is isomorphic to $\Soc N_\pp(D)$ for some simple finite-dimensional $\mm$-module $D$, and $N$ has finite type over $\tt$.
\end{lemma}
\begin{proof}
Let $\{v_1,...,v_n\}$ generate $N$ over $U(\gg)$ and let $B=(U(\bar{\pp}))v_1+(U(\bar{\pp}))v_2+...+(U(\bar{\pp}))v_n$. Then $B$ is a finite-dimensional $\pp$-module. Moreover, $N$ is a quotient of $\ind_{\bar{\pp}}^\gg B$, for which the lemma is well known.
\end{proof}
\begin{prop}\label{pr2.9}
Let $N\in\mathcal{C}_{\bar{\pp},\tt}$. Then $R^{i}\Gamma_{\kk,\tt}(N)$ has finite length for all $i$.
\end{prop}
\begin{proof}
By \refle{le2.8}, $N_{\tt}^*\in\mathcal{C}_{\pp,\tt}$. Any module in $\mathcal{C}_{\pp,\tt}$ admits a resolution by modules of the form $\ind_\pp^\gg C.$, where each $C_{k}$ is a finite-dimensional $(\pp,\tt)$-module: $\ind_\pp^\gg C.\rightarrow N^*_{\tt}\rightarrow 0$. By considering the $\tt$-finite dual of this resolution, we obtain a resolution of $(N^*_{\tt})^*_{\tt}\cong N$ by modules of the form $N_\pp(A^{\cdot})$, where each $A^k$ is a finite-dimensional $(\pp,\tt)$-module.

Write this resolution as $N\hookrightarrow N_\pp(A^{\cdot})$. We have a convergent spectral sequence of $\gg$-modules with $E_2^{a,b}=R^a\Gamma_{\kk,\tt}(N_\pp(A^b))$, and which abuts to $R^{\cdot}\Gamma_{\kk,\tt}(N)$. If we fix $i\in\field{Z}_{\geq 0}$, there are only finitely many terms $E_2^{a,b}$ with $a+b=i$, since both $a\geq 0$ and $b\geq 0$. Hence, $\bigoplus_{a+b=i}E_{\infty}^{a,b}$ has finite length. Finally, $R^{i}\Gamma_{\kk,\tt}(N)$ has finite length
\end{proof}

By $\mathcal{C}_{\kk}$ we denote the full subcategory of $\gg$-mod consisting of $(\gg,\kk)$-modules which have finite type over $\kk$ and have finite length over $\gg$.
\begin{theo}
If $N\in\mathcal{C}_{\bar{\pp},\tt}$ and $i\geq 0$, then $R^{i}\Gamma_{\kk,\tt}(N)\in\mathcal{C}_{\kk}$.
\end{theo}
\begin{proof}
The statement follows from \refprop{pr2.9} and from the "finiteness statement" of Theorem 2.4 c) in \cite{Z}.
\end{proof}

If $\mathcal{A}$ is a full abelian subcategory of $\gg$-mod, let $K_0(\mathcal{A})$ be the Grothendieck group of $\mathcal{A}$.
\begin{defi}
If $N\in\mathcal{C}_{\bar{\pp},\tt}$, let $\Theta_{\kk,\tt}(N)=\sum (-1)^i[R^{i}\Gamma_{\kk,\tt}(N)]$ in $K_0(\mathcal{C}_{\kk})$.
\end{defi}
The fact that $\Theta_{\kk,\tt}(N)$ is well-defined follows from the vanishing statement of Theorem 2.4 b) in \cite{Z}.
\begin{prop}
The map $N\mapsto\Theta_{\kk,\tt}(N)$ yields a non-zero homomorphism $\Theta_{\kk,\tt}: K_0(\mathcal{C}_{\bar{\pp},\tt})\rightarrow K_0(\mathcal{C}_{\kk})$.
\end{prop}
\begin{proof}
This is a well-known fact which follows from the long exact sequence of cohomology.
\end{proof}

\textbf{Example}.

a) If $E$ is a finite-dimensional simple $\mm$-module with highest weight $\nu$ such that $\nu+\tilde{\rho}$ is regular and $\gg$-dominant, then $\Theta_{\kk,\tt}(N_\pp(E))=(-1)^s[F^s(\pp,E)]$. If $\mu$ is $\bb_\kk$-dominant and $\kk$-integral, then $\Theta_{\kk,\tt}(N_\pp(E))\neq 0$ by the remark after \refprop{pr2.5}.

b) $\Theta_{\kk,\tt}(\field{C})=|W_{\kk}|[\field{C}]$, where $W_{\kk}$ is the Weyl group of $\kk$. Indeed, in the proof of Theorem 2.4 in \cite{Z} it is shown that $$\Hom_\kk(V,R^i\Gamma_{\kk,\tt}(\field{C}))\cong\Ext^i_{\kk,\tt}(V,\field{C})$$ for any simple finite-dimensional $\kk$-module $V$. Since $\Ext^i_{\kk,\tt}(V,\field{C})=0$ for $V\not\simeq\field{C}$, we conclude that $\Theta_{\kk,\tt}(\field{C})=\sum_{i}(-1)^i\dim\Ext^i_{\kk,\tt}(\field{C},\field{C})$. Moreover, $\Ext^i_{\kk,\tt}(\field{C},\field{C})=H^i(\kk,\tt,\field{C})$, where $H^i(\kk,\tt,\field{C})$ stands for the relative Lie algebra cohomology. 

It is well-known that $H^i(\kk,\tt,\field{C})$ is the cohomology of the variety $K_0/T_0$, $K_0$ being a connected affine real algebraic group with Lie algebra $\kk$ and $T_0$ being a torus in $K_0$ with Lie $T_0=\tt_0$. Moreover, $K_0/T_0$ is homeomorphic to the flag variety of $K$ and hence the Euler characteristic of $K_0/T_0$ is $|W_\kk|$, by the Bruhat decomposition of the flag variety. Thus, $$\sum_i(-1)^i\dim H^i(\kk,\tt,\field{C})=|W_\kk|.$$

\section{On the $\nn$-cohomology of $(\gg,\kk)$-modules}

We start by recalling [PZ2, Proposition 1]  in the case when $\kk$ is semisimple. In what follows, $M$ will denote a $(\gg,\kk)$-module. Note that $M$ is automatically in the class $\mathcal{M}$.
\begin{prop}
In the category of $\tt$-weight modules, there exists a bounded (not necessarily first quadrant) cohomology spectral sequence which converges to $H^\cdot(\nn,M)$, with 
\[
E_1^{a,b}=H^{a+b-R(a)}(\nn_\kk,M)\otimes V_a^*,
\]
where $a$ runs over $\{0,\dots,y\}$ for some $y, R$ is a monotonic function on $\{0,\dots,y\}$ with values in $\ZZ_{>0}$ such that $R(a)\leq a$ and $R(y)=r$, $V_a$ is a $\tt$-submodule of $\Lambda^{R(a)}(\nn\cap\kk^\perp)$ for every $a$, and $V_y=\Lambda^r(\nn\cap\kk^\perp)$. We also have $\bigoplus_{R(a)=p}V_a=\Lambda^p(\nn\cap\kk^\perp)$. 
\end{prop}
Suppose we are interested in $H^j(\nn,M)$ for a fixed $j$. Write 
\[
E_1^j:=\bigoplus_{p=0}^r H^{j-p}(\nn_\kk,M)\otimes\Lambda^p(\nn\cap\kk^\perp)^*.
\]
Then $E_1^j=\bigoplus_{a+b=l}E_1^{a,b}$.

\begin{lemma} \label{le1}
Fix $\varkappa \in \tt^*$ and $j$, $0\leq j\leq\dim\nn=n$. Assume that $(E_1^{j-1})^{\varkappa}=(E_1^{j+1})^{\varkappa}=0$. Then 
\begin{equation}\label{eqn3}
H^j(\nn,M)^\varkappa\simeq (E_1^j)^\varkappa=\bigoplus_{p=0}^r\left(H^{j-p}(\nn_\kk,M)\otimes \Lambda^p(\nn\cap\kk^\perp)^*\right)^\varkappa.
\end{equation}
\end{lemma}
\begin{proof}
This follows immediately from the definition of a convergent spectral sequence of vector spaces. 
\end{proof}

As a special case we have the following lemma. Recall that $s:=\dim(\nn\cap\kk)$, $r:=\dim(\nn\cap\kk^\perp)$.
\begin{lemma}\label{le3.5}
If $H^s(\nn_\kk,M)^{\varkappa'}=0$ for $\varkappa':=\varkappa+2\rho_\nn^\perp$, then $H^n(\nn,M)^\varkappa=0$.
\end{lemma}
\begin{proof}
The isomorphism \refeq{eqn3} implies $H^n(\nn,M)^\varkappa\simeq H^s(\nn_\kk,M)^{\varkappa '}\otimes\Lambda^r((\nn\cap\kk^\perp)^*)$.
\end{proof}

As a consequence we have the following.
\begin{prop}
If $\varkappa '=\varkappa+2\rho_\nn^\perp$ is $\kk$-dominant integral, then $H^n(\nn,M)^\varkappa=0$.
\end{prop}
\begin{proof}
Kostant's theorem \cite{Ko} implies that if $\eta$ is $\kk$-dominant integral, $H^s(\nn_\kk,V(\eta))$ has pure weight $-w_\mm(\eta)-2\rho$, where $w_\mm$ is the longest element of $W_\kk$. We have as a consequence that $H^s(\nn_\kk,V(\eta))^{\varkappa'}=0$ if $\varkappa'$ is $\kk$-dominant integral. Since $M$ is a direct sum of simple finite dimensional $\kk$-modules $V(\eta)$ for various $\eta$, we conclude that $H^s(\nn_\kk,M)^{\varkappa'}=0$. (We have used our assumption that $\rk \kk_{ss}>0$). Then $H^\nn(\nn,M)^\varkappa=0$ by \refle{le3.5}.
\end{proof}

We now recall that $H^\cdot(\nn,M)$ is an $(\mm,\mm\cap\kk)$-module. This is established in [V2, Ch. 5] in the case when $\kk$ is a symmetric subalgebra but the argument extends to the case of a general reductive in $\gg$ subalgebra $\kk$. Note that $\mm\cap\kk=\tt$. The following statement is identical to [PZ2, Corollary 3].
\begin{prop}\label{pr3.5}
a) If $M$ is a $(\gg,\kk)$-module of finite type, then $H^\cdot(\nn,M)$ is an $(\mm,\tt)$-module of finite type. Moreover, if $M$ is $Z_{U(\gg)}$-finite (i.e. the action of $Z_{U(\gg)}$ on $M$ factors through a finite-dimensional quotient of  $Z_{U(\gg)}$ ) then $H^\cdot(\nn,M)$ is  $Z_{U(\mm)}$-finite.

b) If $\pp$ is a minimal compatible parabolic subalgebra and $M$ is a $(\gg,\kk)$-module of finite type which is in addition  $Z_{U(\gg)}$-finite, then $H^\cdot(\nn,M)$ is finite dimensional.
\end{prop}

\section{Reconstruction of $(\gg,\kk)$-modules}

Suppose $M$ is simple. Let $V(\mu)$ be a minimal $\kk$-type of $M$ (a priori $V(\mu)$ is not unique). Suppose $\mu+2\rho$ is $(\gg,\kk)$-regular and $\pp=\pp_{\mu+2\rho}$.
\begin{defi}
The pair $(M,\mu)$ as above is \emph{strongly reconstructible} if $H^r(\nn,M)^{\mu-2\rho_\nn^\perp}$ is a simple $\mm$-module and there is an isomorphism of $\gg$-modules
\begin{equation}\label{eqn4}
M\simeq\Soc F^s(\pp,H^r(\nn,M)^{\mu-2\rho_\nn^\perp}).
\end{equation}
\end{defi}
The isomorphism \refeq{eqn4} implies via Theorem 2 of \cite{PZ3} that $V(\mu)$ is the unique minimal $\kk$-type of M. Moreover, $\dim\Hom_\kk(V(\mu),M)\leq\dim E$. Therefore, if the pair $(M,\mu)$ is strongly reconstructible, $\mu$ is determined by $M$ as the highest weight of the unique minimal $\kk$-type of $M$. This allows us to simply speak of strongly reconstructible simple $(\gg,\kk)$-modules rather than of strongly reconstructible pairs.

The first reconstruction theorem of \cite{PZ2} now implies the following.
\begin{theo}
If $M$ is a simple $(\gg,\kk)$-module of finite type with a generic minimal $\kk$-type $V(\mu)$, then $M$ is strongly reconstructible. 
\end{theo}

Below we will see (in particular in Subsection 8.2) that the converse to the above theorem is false. We will also see (in Subsection 8.7) examples of simple finite-dimensional $\mm$-modules $E$ such that $F^s(\pp,E)$ has a reducible socle. 
\begin{defi}
A simple $(\gg,\kk)$-module $M$ of finite type over $\kk$ is weakly reconstructible if for some minimal $\kk$-compatible parabolic subalgebra $\pp$, there exists an injective homomorphism of $\gg$-modules 
\[
M\hookrightarrow F^s(\pp,\Top H^r(\nn,M)).
\]
\end{defi}
\begin{theo}
Let $M_1,M_2\in\mathcal{C}_{\kk}$ be simple with central characters $\theta_{\lambda_1}$ and $\theta_{\lambda_2}$ respectively. Assume $\lambda_1$ and $\lambda_2$ are dominant regular with respect to a Borel subalgebra $\bb\subset\gg$; assume further that $\lambda_{2}-\lambda_{1}$ is dominant integral . Finally assume that $\Phi$ is a translation functor such that $\Phi(M_1)\cong M_2$. Then, $M_2$ is weakly reconstructible if and only if $M_1$ is weakly reconstructible.
\end{theo}
\begin{proof}
Assume that $M_1$ is weakly reconstructible. Then for some simple quotient $E_1$ of $H^r(\nn,M_1)$, we have an injection of $\gg$-modules $\alpha_1 :M_1\rightarrow F^s(\pp,E_1)$. By assumption, $\Phi$ is an equivalence of categories. Hence we have an injection $\alpha_2 :M_2\rightarrow \Phi (F^s(\pp,E_1))$. By \cite{Z}, we have an isomorphism $\Phi (F^s(\pp,E_1))\cong F^s(\pp,E_2)$ for a simple $\mm$-module $E_2$. Moreover, we have a translation functor $\Phi^\mm$ such that $\Phi^\mm (E_1)\cong E_2$. Thus, we have an injection $\alpha_2^1 :M_2\rightarrow F^s(\pp,\Phi^\mm (E_1))$.

Now let $\theta^\mm_{\chi_i}$ be the central character of  $E_i$ for $i=1,\; 2$. Let $P_{\chi_1}^\mm$ and $P_{\chi_2}^\mm$ be the respective projection functors. By assumption we have a surjection  of $\mm$-modules $\beta_1:P_{\chi_1}^\mm (H^r(\nn,M_1))\rightarrow E_1$. Apply the translation functor $\Phi^\mm$ to $\beta_1$ to obtain a surjection $\beta_2:\Phi^\mm P_{\chi_1}^\mm (H^r(\nn,M_1))\rightarrow \Phi^\mm (E_1)\cong E_2$. By [KV, Ch. 7], $\Phi^\mm P_{\chi_1}^\mm (H^r(\nn,M_1))\cong P_{\chi_2}^\mm (H^r(\nn,M_1))$. Thus, $E_2$ is a quotient of $H^r(\nn,M_2)$. Finally, the injection $\alpha_2^1$ yields an injection $M_2\rightarrow F^s(\pp, \Top H^r(\nn,M_2))$. Hence $M_2$ is weakly reconstructible.
\end{proof}

\section{Preliminary results on $(\gg,sl(2))$-modules}
From now on we assume that $\kk$ is isomorphic to $sl(2,\CC)$. We fix a standard basis $\{e,h,f\}$ for $\kk$; the eigenvalues of $\ad h$ in $\gg$ will be integers. Let $\tt=\CC h$ be the Cartan subalgebra of $\kk$ generated by $h$, and let $\pp=\pp_{h^*}$ where $h^*\in\tt^*$, $h^*(h)=1$. The subalgebra $\pp$  is automatically a minimal $\tt$-compatible parabolic subalgebra.

For $\kk\simeq sl(2)$, our \refle{le1} simplifies to the following.
\begin{lemma}\label{le2}
Fix $\varkappa\in\tt^*$ and $j$, $0\leq j\leq r+1$. Write $$E_1^j=H^0(\nn_\kk,M)\otimes \Lambda^j(\nn\cap\kk^\perp)^*\oplus H^1(\nn_\kk,M)\otimes \Lambda^{j-1}(\nn\cap\kk^\perp)^*.$$ Suppose that $(E_1^{j-1})^{\varkappa}=(E_1^{j+1})^{\varkappa}=0$. Then there is an isomorphism of $\tt$-modules $$H^j(\nn,M)^\varkappa\simeq (E_1^{j})^{\varkappa}\cong (H^0(\nn_\kk,M)\otimes \Lambda^j(\nn\cap\kk^\perp)^*\oplus H^1(\nn_\kk,M)\otimes \Lambda^{j-1}(\nn\cap\kk^\perp)^*)^{\varkappa}.$$
\end{lemma}

In particular, let $j=r$ and assume $\varkappa '=\varkappa+2\rho_\nn^\perp$ is dominant integral for $\kk$. Then $(E_1^{r-1})^{\varkappa}=0$ implies 
\begin{equation}\label{eqn5}
H^r(\nn,M)^\varkappa\simeq H^0(\nn_\kk,M)^{\varkappa '}\oplus \left(H^1(\nn_\kk,M)\otimes (\nn\cap\kk^\perp)\right)^{\varkappa '}.
\end{equation}
More generally, $\varkappa '$ dominant integral implies 
\begin{equation}
\dim H^r(\nn,M)^\varkappa\leq \dim H^0(\nn_\kk,M)^{\varkappa '}+\dim \left( H^1(\nn_\kk,M)\otimes (\nn\cap\kk^\perp)\right)^{\varkappa '}. 
\end{equation}  

From now on we identify integral weights $\varkappa$ of $\tt$ with the corresponding integers, $\varkappa(h)$. Let $\lambda_1$ and $\lambda_2$ be the maximum and submaximum weights of $\tt$ in $\nn\cap\kk^\perp$ (we consider $\lambda_1$ and $\lambda_2$ as integers); if $\lambda_1$ has multiplicity at least two in $\nn\cap\kk^\perp$, then $\lambda_2=\lambda_1$.
\begin{prop}\label{pr5.2}
Let $\mu$ be a nonnegative integer and let $M$ be a $(\gg,\kk)$-module with the property that $\delta<\mu$ implies $M[\delta]=0$.

\noindent a) If $\mu\geq \half \lambda_1$, then $\dim H^r(\nn,M)^\omega\leq \dim H^0(\nn_\kk, M)^\mu.$

\noindent b) If $\mu\geq \half(\lambda_1+\lambda_2)$, then 
\[
\dim H^r(\nn,M)^\omega = \dim H^0(\nn_\kk,M)^\mu.
\]
\end{prop}
\begin{proof}
\noindent a)
Our hypothesis on $M$ implies that if $M[\delta]\neq 0$, then $\delta\geq\mu\geq\half\lambda_1$. Since $H^1(\nn_\kk,M)$ has weights $-\delta-2$ with $\delta$ as above, we see that \newline $\left( H^1(\nn_\kk,M)\otimes(\nn\cap\kk^\perp)\right)^\mu=0$. Hence (2) implies the inequality in a).

\noindent b) It suffices to show that $(E_1^{r-1})^\omega=0$. Then the statement from \refeq{eqn5}, taking into account the vanishing of $(H^1(\nn_\kk,M)\otimes(\nn\cap\kk^{\perp}))^\mu$, implies (b).

We now check that $(E_1^{r-1})^\omega=0$. We have $$(E_1^{r-1})^\omega = (H^0(\nn_\kk,M)\otimes\Lambda^{r-1}((\nn\cap\kk^{\perp})^*))^\omega\oplus(H^1(\nn_\kk,M)\otimes\Lambda^{r-2}((\nn\cap\kk^\perp)^*))^\omega.$$
Furthermore,
$$\Lambda^{r-1}((\nn\cap\kk^\perp)^*)\cong(\nn\cap\kk^\perp)\otimes\Lambda^r((\nn\cap\kk^\perp)^*),$$
$$\Lambda^{r-2}((\nn\cap\kk)^*)\cong\Lambda^2(\nn\cap\kk^\perp)\otimes\Lambda^r((\nn\cap\kk^\perp)^*)$$
implies
\begin{equation}\label{eqn6}
(E_1^{r-1})^\omega = (H^0(\nn_\kk,M)\otimes(\nn\cap\kk^\perp))^\mu\oplus(H^1(\nn_\kk,M)\otimes\Lambda^2(\nn\cap\kk^\perp))^\mu
\end{equation}
as the weight of $\Lambda^r(\nn\cap\kk^\perp)$ equals $2\rho_\nn^\perp$. The first term of \refeq{eqn6} vanishes as the $\tt$-weights of $\mm\cap\kk^\perp$ are strictly positive and the smallest $\tt$-weight of $H^0(\nn_\kk,M)$ is $\mu$. The maximal $\tt$-weight of the second term is $-\mu -2+\lambda_1+\lambda_2$, hence the inequality $\mu\geq\frac{\lambda_1+\lambda_2}{2}$ implies the vanishing of the second term of \refeq{eqn6}.
\end{proof}

Our next task is to state and prove a vanishing theorem for $F^0(\pp,E)$, where $E$ is a simple finite dimensional $\mm$-module. Let $\omega\in\tt^*$ be the weight of $\tt$ in $E$. 
\begin{prop}\label{pr5.3}
Suppose $\mu=\omega+2\rho_\nn^\perp$ and $\mu\geq 0$. Then $F^0(\pp,E)=0$.
\end{prop}
\begin{proof}
By definition 
\[
F^0(\pp,E)=\Gamma_{\kk,\tt}(N_\pp(E)).
\]
We have $N_\pp(E)_\tt^*\cong\ind_\pp^\gg(E^*\otimes\Lambda^{\dim(\nn)}(\nn^*))$ and the $\bb$-highest weight of $\ind_\pp^\gg(E^*\otimes\Lambda^{\dim\nn}(\nn^*))$ equals $-\nu -2\tilde{\rho}_\nn\in\hh^*$.

On the other hand, $\nu '(h)\geq 0$ for any $\bb$-dominant weight. This follows from the fact that any $\bb$-dominant weight is a non-negative linear combination of roots of $\bb$ (see for instance \cite{Kn}, p. 686).

The $\gg$-module $N_\pp(E)$ has a finite-dimensional submodule if and only if $N_\pp(E)_\tt^*$ has a finite-dimensional quotient. Note that $\Gamma_{\kk,\tt}(N_\pp(E))$ is an integrable $\gg$-module as $N_\pp(E)$ is $\bar{\pp}$-locally finite and $\kk$ and $\bar{\pp}$ generate $\gg$. Therefore, $\Gamma_{\kk,\tt}(N_\pp(E))=0$ whenever $N_\pp(E)$ has no finite-dimensional submodule, i.e. whenever $-\gamma-2\tilde{\rho}_\nn$ is not $\bb$-dominant. The fact that $(-\gamma-2\tilde{\rho}_\nn)(h)=-\omega-2\rho_\nn=\mu-2\rho <0$, allows us to conclude that $\nu-2\tilde{\rho}_\nn$ is not $\bb$-dominant, i.e. that $\Gamma_{\kk,\tt}(N_\pp(E))=0$.
\end{proof}
\begin{prop}\label{pr5.4}
$F^2(\pp,E)=0$. 
\end{prop}
\begin{proof}
The statement is a direct corollary of Proposition 3 a) in \cite{PZ2}. Note that this proof does not use genericity.
\end{proof}
\begin{prop}\label{pr5.5}(cohomological Frobenius reciprocity) Let $\mu\geq 0$. If $M$ is a $(\gg,\kk)$-module such that $H^\cdot(\nn,M)$ is finite dimensional, then we have a natural isomorphism $$\Hom_\gg(M,F^1(\pp,E))\cong\Hom_\mm(H^r(\nn,M)^\omega,E).$$
\end{prop}
\begin{proof}
This follows from the existence of a (not necessarily first quadrant) spectral sequence with $E_2$-term 
$E_2^{a,b}=\Ext_{\mm,\tt}^a(H^{r-b}(\nn,M),E)$ converging to $\Ext_{\gg,\kk}^{a+b}(M,F^1(\pp,E))$, see Proposition 6 of \cite{PZ2}. By assumption $H^\cdot(\nn,M)$  is finite dimensional. Choose $b_0$ to be the least possible integer with $\Ext^\cdot_{\mm,\tt}(H^{r-b_0}(\nn,M),E)\neq 0$. 

By the same argument as in the proof of Theorem 2, b) in \cite{PZ2}, we conclude that $\Hom_\mm(H^{r-b_0}(\nn,M),E)\neq 0$. Thus, $E_2^{0,b_0}\neq 0$ and $E_2^{a,b}=0$ for $b<b_0$. Consequently, $E_2^{0,b_0}\cong E_\infty^{0,b_0}$ and we deduce that $\Ext_{\gg,\kk}^{b_0}(M,F^1(\pp,E))\neq 0$. Hence, $b_0\geq 0$ and the spectral sequence is a  first quadrant spectral sequence, with corner isomorphism $\Hom_\gg(M,F^1(\pp,E))\cong\Hom_\mm(H^r(\nn,M),E)$.
\end{proof}
\begin{corollary}\label{cor5.6}
Suppose $\mu\in\field{Z}_{\geq 0}$.

a) Let $X$ be any $\gg$-submodule of $F^1(\pp,E)$. Then $E$ is a quotient of $H^r(\nn,X)^\omega$. In particular, if $X$ is simple, then $X$ is weakly reconstructible.

b) Let $M$ be a simple $(\gg,\kk)$-module such that $H^\cdot(\nn,M)$ is finite dimensional and $E$ is isomorphic to a quotient of $H^r(\nn,M)^\omega$. Then $M$ is isomorphic to a submodule of $\Soc F^1(\pp,E)$. In particular $M$ is weakly reconstructible and $M$ has finite type over $\kk$.
\end{corollary}
\begin{corollary}
Fix a central character $\theta$. The set of isomorphism classes of simple $(\gg,\kk)$-modules $M$ with central character $\theta$ such that $\dim H^r(\nn,M)<\infty$ and $H^r(\nn,M)^{\varkappa}\neq 0$ for some $\varkappa\in\field{Z}_{\geq -2}$ is finite.
\end{corollary}
\begin{proof}
A $\gg$-module $M$ as in the corollary is isomorphic by \refcor{cor5.6} (b) to a $\gg$-submodule of $F^1(\pp,E')$, where $E'$ runs over finitely many simple finite-dimensional $\pp$-modules. Since $F^1(\pp,E')$ has finite length for each $E'$ by \refth{th2.6}, the statement follows.
\end{proof}
\begin{theo}\label{th5.8}
Suppose $\kk$ is regular in $\gg$. Let $M$ be a simple $(\gg,\kk)$-module, not necessarily of finite type over $\kk$, with lowest $\kk$-type $V(\mu)$ for $\mu\geq 0$.

a) If $H^r(\nn,M)^\omega\neq 0$, there exists a 1-dimensional simple quotient $E$ of $H^r(\nn,M)^\omega$. For any such $E$ we have an injective homomorphism $M\rightarrow F^1(\pp,E)$. Hence $M$ is weakly reconstructible and $M$ is of finite type over $\kk$.

b) If $M$ is of infinite type over $\kk$, then $H^r(\nn,M)^\omega =0$.\footnote{See Theorem 9 of \cite{PZ1}.}
\end{theo}
\begin{proof}
a) Since $M$ is simple, \refprop{pr3.5} b) implies that $H^\cdot(\nn,M)$ is finite dimensional. Note that the regularity of $\kk$ in $\gg$ implies that $\mm$ is a Cartan subalgebra of $\gg$. Hence there exists a 1-dimensional simple $\mm$-quotient $E$ of $H^r(\nn,M)^\omega$. \refprop{pr5.5} implies now that any such $E$ induces an injective homomorphism of $M$ into $F^1(\pp,\Top H^r(\nn,M))$. In particular, $M$ is weakly reconstructible and is of finite type over $\kk$.

b) Follows from a).
\end{proof}

\section{Strong reconstruction of $(\gg,sl(2))$-modules}

In the rest of the paper $E,\pp,\mm,\tt,\mu,\omega$ are as in Section 1. We start with the following result on the  of $(\gg,\kk)$-modules. 
\begin{theo}\label{th6.1}
Let $\mu\geq\frac{1}{2}\lambda_1$. Then $\Soc F^1(\pp,E)=\bar{F}^1(\pp,E)$ and $\bar{F}^1(\pp,E)$ is simple. In particular, $\dim\Hom_\kk(V(\mu),\Soc F^1(\pp,E))=\dim E$.
\end{theo}
\begin{proof}
Let $X$ be a non-zero submodule of $F^1(\pp,E)$. Since $F^1(\pp,E)$ is a $(\gg,\kk)$-module of finite type and is $Z_{U\gg}$-finite, \refprop{pr3.5} b) and \refprop{pr5.5} apply to $X$, yielding a surjective homomorphism of $\mm$-modules $H^r(\nn,X)^\omega\rightarrow E$. Hence $\dim E\leq \dim H^r(\nn,X)^\omega$. Next, by \cite{PZ3}, $$\dim H^0(\nn_\kk,F^1(\pp,E))^\omega=\dim E,$$ therefore $\dim H^0(\nn_\kk,X)^\omega\leq\dim E$ by the left exactness of $H^0(\nn_\kk,\cdot)$. Finally, by \refprop{pr5.2} a), $\dim H^r(\nn,X)^\omega\leq\dim H^0(\nn_\kk,X)^\mu$. Combining these inequalities we see that $$\dim H^r(\nn,X)^\omega =\dim H^0(\nn_\kk,X)^\mu =\dim E.$$ Hence $X[ \mu ]=F^1(\pp,E)[ \mu ]$, or equivalently $X\supseteq \bar{F}^1(\pp,E)$.
Since in this way $\bar{F}^1(\pp,E)$ is contained in any non-zero submodule of $F^1(\pp,E)$, $\bar{F}^1(\pp,E)$ is simple and $\bar{F}^1(\pp,E)=\Soc F^1(\pp,E)$.
\end{proof}
\begin{corollary}\label{cor6.2}
Under the assumptions of \refth{th6.1}, let $X\neq 0$ be a $\gg$-submodule of $F^1(\pp,E)$. Then the minimal $\kk$-type of $X$ is $V(\mu)$, $\dim\Hom_\kk(V(\mu),X)=\dim E$, and there is an isomorphism of $\mm$-modules $H^r(\nn,X)^\omega\cong E$.
\end{corollary}
\begin{proof}
The statement was established in the proof of \refth{th6.1}.
\end{proof}
\begin{corollary}\label{cor6.3}
Let $M$ be a simple $(\gg,\kk)$-module whose minimal $\kk$-type $V(\mu)$ satisfies $\mu\geq\frac{1}{2}\lambda_1$. Then, if $H^\cdot(\nn,M)$ is finite-dimensional and $H^r(\nn,M)^\omega\neq 0$, $M$ is strongly reconstructible. 
\end{corollary}
\begin{proof}
Let $E'$ be a simple quotient of the $\mm$-module $H^r(\nn,M)^\omega$. By \refprop{pr5.5}, $M$ is a simple submodule of $F^1(\pp,E')$, hence by \refth{th6.1} $$M\cong\Soc F^1(\pp,E').$$
\end{proof}

\begin{corollary}\label{cor6.4}
Let $M$ be a simple $(\gg,\kk)$-module of finite type such that its minimal $\kk$-type $V(\mu)$ satisfies $\mu\geq\frac{1}{2}(\lambda_1 + \lambda_2)$. Then $H^\cdot(\nn,M)$ is finite-dimensional, and $H^r(\nn,M)^\omega\neq 0$, hence $M$ is strongly reconstructible by \refcor{cor6.3}.
\end{corollary}
\begin{proof}
The statement follows from \refcor{cor6.4} via \refprop{pr5.2} b).
\end{proof}

\begin{corollary}\label{cor6.5}
The correspondences $$M\rightsquigarrow H^r(\nn,M)^\omega$$
$$E\rightsquigarrow\Soc F^1(\pp,E)$$
induce mutually inverse bijections between the set of isomorphism classes of simple $(\gg,\kk)$-modules of finite type $M$ whose minimal $\kk$-type $V(\mu)$ satisfies $\mu\geq\frac{1}{2}(\lambda_1 + \lambda_2)$ and the set of isomorphism classes of finite dimensional $\mm$-modules on which $\tt$ acts via $\omega = \mu - 2\rho_\nn^\perp$, where $\field{Z}_{\geq 0}\ni\mu\geq\frac{1}{2}(\lambda_1 + \lambda_2)$.
\end{corollary}
\begin{corollary}\label{cor6.6}
Suppose $\kk$ is regular in $\gg$ (i.e. let $\tt$ contain an element regular in $\gg$). Suppose $M$ is a simple $(\gg,\kk)$-module (not necessarily of finite type over $\kk$) with lowest $\kk$-type $V(\mu)$. If $\mu\geq\frac{1}{2}(\lambda_1 + \lambda_2)$, then $E=H^r(\nn,M)^\omega$ is a 1-dimensional $\mm$-module and $M\cong \bar{F}^1(\pp,E)$. Thus $M$ is strongly reconstructible. In particular, $M$ has finite type over $\kk$.
\end{corollary}
\begin{proof}
We apply \refth{th5.8} a) and \refprop{pr5.2} b) to conclude that, for any 1-dimensional quotient $E$ of $H^r(\nn,M)^\omega$, we have an injection $M\rightarrow F^1(\pp,E)$. Since $\mu\geq\frac{1}{2}(\lambda_1 + \lambda_2)$, there are isomorphisms $M\cong\bar{F}^1(\pp,E)$ and $H^r(\nn,M)^\omega\cong E$ (the latter is an isomorphism of $\mm$-modules).
\end{proof}

\textbf{Example.} Let $\gg$ be a classical simple Lie algebra of rank $n$ and $\kk\cong sl(2)$ be a principal subalgebra. Then the claim of \refcor{cor6.5} is proved in \cite{PZ2} under the assumption that $\mu+1\geq 2(\sum_i r_i)$, where $\tilde{\rho}=\frac{1}{2}\sum_i r_i\alpha_i$, $\alpha_i$ being the simple roots of $\bb$. It is well-known that $2(\sum_i r_i)$ grows cubically with the growth of $n$, while the value $\frac{1}{2}(\lambda_1+\lambda_2)$ has linear growth in $n$. Therefore, for large $n$, the result of \refcor{cor6.5} strengthens considerably Theorem 3 of \cite{PZ2} for $\kk$ being a principal $sl(2)$-subalgebra . On the other hand, we will see in Section 8 that for $n=2$, there are cases where the bound $2(\sum_i r_i)-1$ is lower that $\frac{1}{2}(\lambda_1+\lambda_2)$.

Set now $\tilde{\kk}:=\kk\oplus C(\kk)$ and note that $\tilde{\kk}$ is a reductive in $\gg$ subalgebra. Recall that $C(\kk)_{ss}\subset \mm_{ss}$. Moreover, $\mm_{ss}\subset C(\kk)\Leftrightarrow \mm_{ss}=C(\kk)_{ss}$.
\begin{prop}
If $\mm_{ss}=C(\kk)_{ss}$, then for any simple $(\gg,\tilde{\kk})$-module $M$ of finite type over $\tilde{\kk}$, $H^\cdot(\nn,M)$ is finite dimensional.
\end{prop}
\begin{proof}
By \refprop{pr3.5} a), $H^\cdot(\nn,M)$ is an $(\mm,\mm\cap\tilde{\kk})$-module of finite type as $\pp$ is $\tilde{\kk}$-compatible (see also [V2, Corollary 5.2.4]). But $\mm\cap\tilde{\kk}=(Z_\mm\cap\tilde{\kk})\oplus\mm_{ss}$. Hence $H^\cdot(\nn,M)$ is an integrable $\mm$-module. Finally, [PZ2, Corollary 3 a)] implies now that $H^\cdot(\nn,M)$ is finite dimensional.
\end{proof}

We conclude this section with some applications to the case when $\tilde{\kk}$ is a symmetric subalgebra.
\begin{prop}\label{pr6.8}
Assume that $\gg$ is simple and $\tilde{\kk}$ is symmetric.

a) If $\gg$ is classical with $rank\geq 4$, the only case of a symmetric pair of the form $(\gg,\tilde{\kk})$ for which $\mm_{ss}$ is not equal to $C(\kk)_{ss}$ is the series $(so(2n), so(3)\oplus so(2n-3))$, where $\kk=so(3)$.

b) If $\gg$ is exceptional, then $\mm_{ss}=C(\kk)_{ss}$. In fact, $\tilde{\kk}$ is symmetric if and only if $\kk$ is conjugate to the $sl(2)$-subalgebra of a highest root of $\gg$.
\end{prop}
\begin{proof}
Follows from the classification of symmetric pairs.
\end{proof}
\begin{corollary}
If $\rk\tilde{\kk}=\rk\gg$ and $\tilde{\kk}$ is symmetric, then $\mm_{ss}=C(\kk)_{ss}$.
\end{corollary}
\begin{proof}
Follows from \refprop{pr6.8}, but a more elegant proof is based on Borel - De Siebenthal \cite{BdS}.
\end{proof}

\begin{corollary}
Assume that $\kk$ is symmetric and $\mm_{ss}=C(\kk)_{ss}$. Let $M$ be a simple $(\gg,\tilde{\kk})$-module.

a) If $H^r(\nn,M)^{\omega}\neq 0$, then $M$ is strongly reconstructible as a $(\gg,\kk)$-module; in particular $M$ has finite type over $\kk$.

b) If $\mu\geq\frac{1}{2}(\lambda_1+\lambda_2)$, then $M$ is strongly reconstructible as a $(\gg,\kk)$-module; in particular $M$ has finite type over $\kk$.
\end{corollary}
 
\section{$\kk$-characters and composition multiplicities of the fundamental series of $(\gg,sl(2))$-modules}

Assume $\mu\in\field{Z}$. Set $L_\pp(E) = \Soc N_\pp(E)$ and recall that $L_\pp(E)$ is simple. Also note that $N_\pp(E)$ and $L_\pp(E)$ are objects of $\mathcal{C}_{\bar{\pp},\tt}$. Denote by $D$ a variable simple finite-dimensional $\pp$-module on which $\tt$ acts via $\mu_D-2\rho_\nn^\perp$. Non-negative integers $m(E,D)$ are determined from the equality $[N_\pp(E)] = \sum m(E,D)[L_\pp(D)]$ in the Grothendieck group $K_0(\mathcal{C}_{\bar{\pp},\tt})$.  We arrange the integers $m(E,D)$ into a matrix $(m(E,D))$ with rows indexed by all possible $E$ and columns indexed by all possible $D$; the rows and columns of $(m(E,D))$ are finitary, i.e. each row and each column have finitely many non-zero entries. The algorithm for computing the integers $m(E,D)$ is discussed in \cite{CC}.
\begin{lemma}\label{le7.1}
Suppose $D$ is not isomorphic to $E$. Then $m(E,D)>0$ implies $\mu_D\geq\mu$.
\end{lemma}
\begin{proof}
We claim that the minimum $\tt$-weight of $N_\pp(E)$ is $\mu+2$. To see this it suffices to note that $N_\pp(E)^*_\tt\simeq\ind_\pp^\gg(E\otimes\Lambda^{\dim\nn}(\nn))^*\cong U(\bar{\nn})\otimes(E\otimes\Lambda^{\dim\nn}(\nn))^*$, as the maximum weight of $N_\pp(E)_\tt^*$ is $-\omega -2\rho_\nn=-\omega -2\rho_\nn^\perp +2\rho_\nn^\perp -2\rho_\nn=-\mu -2\rho=-\mu -2$ (note that $\rho=1$). Thus, if $L_\pp(D)$ is a composition factor of $N_\pp(E)$ and $D\not\cong E$, we have $\mu_D+2>\mu+2$, or equivalently, $\mu_D>\mu$.
\end{proof}
\begin{prop}\label{pr7.2} 
Let $\mu\geq 0$. Then $[F^1(\pp,E)] = \sum_D m(E,D) [R^1 \Gamma_{\kk,\tt} (L_\pp(D))]$.
\end{prop}
\begin{proof}
Taking into account that $\Theta$ is a homomorphism, it suffices to prove:

a) $F^i(\pp,E)=R^i\Gamma_{\kk,\tt}(N_\pp(E))=0$ for $i\neq 1$ and $\mu\geq 0$;

b) $R^i\Gamma_{\kk,\tt}(L_\pp(D))=0$ for $i\neq 1,\; m(E,D)>0,\; \mu\geq 0$.

Part a) follows from \refprop{pr5.3} and \refprop{pr5.4}. 

To prove part b), note that $m(E,D)>0$ implies $\mu_D\geq\mu\geq 0$ by \refle{le7.1}. Then $F^0(\pp,D)=0$ implies $R^0\Gamma_{\kk,\tt}(L_\pp(D))=0$ as $R^0\Gamma_{\kk,\tt}(L_\pp(D))\subseteq F^0(\pp,D)$. To see that $R^2\Gamma_{\kk,\tt}(L_\pp(D))=0$, note that, by the Duality Theorem in \cite{EW}, $(R^2\Gamma_{\kk,\tt}(L_\pp(D)))_\kk^*\cong\Gamma_{\kk,\tt}(L_\pp(D)_\tt^*)$. The $\gg$-module $L_\pp(D)_\tt^*$ is $\pp$-locally finite, simple and infinite dimensional. As $\pp$ and $\kk$ generate $\gg$, $\Gamma_{\kk,\tt}(L_\pp(D)_\tt^*)\neq 0$ would imply $\dim L_\pp(D)_\tt^*<\infty$. As the latter is false, $\Gamma_{\kk,\tt}(L_\pp(D)_\tt^*)=0$. Part b) is proved.
\end{proof}
\begin{prop}\label{pr7.3}
Suppose $\mu\geq 0$.

a) Then $R^1\Gamma_{\kk,\tt}(L_\pp(E))\neq 0$, and the lowest $\kk$-type of $R^1\Gamma_{\kk,\tt}(L_\pp(E))$ is $V(\mu)$ of multiplicity $\dim E$.

b) We have the following inclusions of $(\gg,\kk)$-modules: $$\bar{F}^1(\pp,E)\subseteq R^1\Gamma_{\kk,\tt}(L_\pp(E))\subseteq F^1(\pp,E).$$

c) If $N_\pp(E)$ is reducible then the submodule $R^1\Gamma_{\kk,\tt}(L_\pp(E))$ of $F^1(\pp,E)$ is proper and non-zero, and hence $F^1(\pp,E)$ is reducible.
\end{prop}
\begin{proof}
a) By \refprop{pr7.2}, we have under our hypothesis
\begin{equation}\label{eqn7}
F^1(\pp,E)[\mu]=\sum_D R^1\Gamma_{\kk,\tt}(L_\pp(D))[\mu].
\end{equation}
By Theorem 2 from \cite{PZ3}, $F^1(\pp,E)[\mu]\cong(\dim E)V(\mu)$. 

Next, apply \refprop{pr7.2} to $D$. We see that $[R^1\Gamma_{\kk,\tt}(L_\pp(D))]$ is a summand of $[F^1(\pp,D)]$. Thus, $\dim (R^1\Gamma_{\kk,\tt}(L_\pp(D))[\mu])\leq\dim (R^1\Gamma_{\kk,\tt}(N_\pp(D))[\mu])$.

Assume now that $m(E,D)>0$ and $E\not\cong D$. Then $\mu_D>\mu$ by \refle{le7.1}. Applying Theorem 2 from \cite{PZ3} a second time, we have $F^1(\pp,D)[\mu]=0$. 

We conclude that if $m(E,D)>0$ and $D\not\cong E$, $R^1\Gamma_{\kk,\tt}(L_\pp(D))[\mu]=0$. So, from formula \refeq{eqn7} above, we deduce that
\begin{equation}\label{eqn8}
F^1(\pp,E)[\mu]=R^1\Gamma_{\kk,\tt}(L_\pp(E))[\mu].
\end{equation}
This proves part a).

b) By the vanishing theorems a) and b) in the proof of \refprop{pr7.2}, the inclusion of $L_\pp(E)$ into $N_\pp(E)$ yields an injection of $R^1\Gamma_{\kk,\tt}(L_\pp(E))$ into $F^1(\pp,E)$; part b) follows immediately from \refeq{eqn8} and the definition of $\bar{F}^1(\pp,E)$.

c) If $N_\pp(E)$ is reducible then for some $D$ with $\mu_D>\mu\geq 0$, $m(E,D)\neq 0$. So, part c) follows from \refprop{pr7.2} and \refprop{pr7.3} a).
\end{proof}

\begin{conjecture}\label{conj7.4}
If $\mu\geq 0$, then $R^1\Gamma_{\kk,\tt}(L_\pp(E))$ is a semisimple $\gg$-module.
\end{conjecture}

If true, this conjecture would imply that all simple constituents of $R^1\Gamma_{\kk,\tt}(L_\pp(E))$ are weakly reconstructible for $\mu\geq 0$. This would follow from \refcor{cor5.6} a). See Subsection 8.7, Examples 1 and 2 for cases when $R^1\Gamma_{\kk,\tt}(L_\pp(E))$ is reducible.
\begin{theo}\label{th7.4}
Assume $\mu\geq\frac{\lambda_1}{2}$. Then $R^1\Gamma_{\kk,\tt}(L_\pp(E))$ is a simple (in particular, non-zero) submodule of $F^1(\pp,E)$ and $R^1\Gamma_{\kk,\tt}(L_\pp(E))=\bar{F}^1(\pp,E)=\Soc F^1(\pp,E)$.
\end{theo}
\begin{proof}
By \refprop{pr7.3}, we have $\bar{F}^1(\pp,E)\subseteq R^1\Gamma_{\kk,\tt}(L_\pp(E))$. By the Duality Theorem, $(R^1\Gamma_{\kk,\tt}(L_\pp(E)))^*_\kk\cong R^1\Gamma_{\kk,\tt}(L_\pp(E)^*_\tt)$.

Let $\bar{F}^1(\pp,E)^\perp$ be the submodule of $R^1\Gamma_{\kk,\tt}(L_\pp(E)^*_\tt)$ consisting of vectors which are orthogonal to $R^1\Gamma_{\kk,\tt}(L_\pp(E))$ via the above duality. By \refprop{pr7.3}, $\bar{F}^1(\pp,E)[\mu]=R^1\Gamma_{\kk,\tt}(L_\pp(E))[\mu]$. Hence, $\tilde{F}^1(\pp,E)^\perp[-\mu]=0$.

By construction, $\bar{F}^1(\pp,E)^\perp$ is a submodule of $F^1(\bar{\pp},E^*)$. It follows from the proof of \refcor{cor6.2} that $\tilde{F}^1(\pp,E)^\perp=0$ and hence $\bar{F}^1(\pp,E)=R^1\Gamma_{\kk,\tt}(L_\pp(E))$.

The statement of \refth{th6.1} implies the remainder of the proof of \refth{th7.4}.
\end{proof}

\begin{corollary}\label{cor7.5}
If $\mu \geq \frac{\lambda_1}{2}$, then:

a) $[F^1(\pp,E)] = \sum m(E,D)[\bar{F}^1(\pp,D)]$.

b) If $N_\pp(E)$ is irreducible, then $\bar{F}^1(\pp,E)=R^1\Gamma_{\kk,\tt}(L_\pp(E))=F^1(\pp,E)$, and $F^1(\pp,E)$ is irreducible. \footnote{\refcor{cor7.5} b) is a strengthening of \refth{th2.1} under the assumption that $\kk\cong sl(2)$.}

c) $\bar{F}^1(\pp,E)[\mu] = F^1(\pp,E)[\mu](\cong \field{C}^{\dim E}V(\mu))$.

d) $[\bar{F}^1(\pp,E)] = \sum p(E,D)[F^1(\pp,D)]$, where (p(E,D)) is the matrix inverse to $(m(E,D))$.

e) $\ch_{\kk} \bar{F}^1(\pp,E) = \sum p(E,D) \ch_{\kk} F^1(\pp,D)$. (See \cite{PZ1} for a formula for $\ch_\kk F^1(\pp,D)$.)

f) $H^r(\nn,R^1\Gamma_{\kk,\tt}(L_\pp(E)))^{\omega}\cong E$; in particular, the $\gg$-module $R^1\Gamma_{\kk,\tt}(L_\pp(E))$ determines $E$ up to isomorphism.
\end{corollary}
\begin{proof}
a) Apply \refprop{pr7.2} and \refth{th7.4}.

b) If $N_\pp(E)$ is irreducible, then $m(E,D)=0$ for $D\not\cong E$ and $m(E,E)=1$. Now apply \refcor{cor7.5} a).

c) Combine formula \refeq{eqn8} with \refth{th7.4}.

d) Follows from a) and the definition of the matrix $(p(E,D))$.

e) Follows from c).

f) Apply \refcor{cor6.2} and \refth{th7.4}.
\end{proof}

Let $n\in\field{Z}$ and let $\mathcal{C}_{\bar{\pp},\tt,n}$ be the full subcategory of $\mathcal{C}_{\bar{\pp},\tt}$ consisting of $(\gg,\tt)$-modules $N$ whose weight spaces $N^\alpha$ satisfy $\alpha\in\field{Z}$ and $\alpha\geq n+2$. Let $\mathcal{C}_{\kk,n}$ be the full subcategory of $\mathcal{C}_{\kk}$ consisting of $(\gg,\kk)$-modules $M$ with minimal $\kk$-type $V(\mu)$ for $\mu\geq n$.

Assume $N$ is a non-zero object in $\mathcal{C}_{\bar{\pp},\tt,2}$.
\begin{lemma}\label{lemma7.7}
$R^i\Gamma_{\kk,\tt}(N)=0$ for $i=0$ and $2$; $R^1\Gamma_{\kk,\tt}(N)\neq 0$.
\end{lemma}
\begin{proof}
$N$ has a finite composition series with simple subquotients $L_\pp(D)$ in $\mathcal{C}_{\bar{\pp},\tt,2}$. We know that $R^i\Gamma_{\kk,\tt}(L_\pp(D))=0$ for $i=0$ and $2$. Therefore our claim follows from the long exact sequence for right derived functors.
\end{proof}
\begin{prop}\label{pr7.8}
The restriction of $R^1\Gamma_{\kk,\tt}(\cdot)$ to the full subcategory $\mathcal{C}_{\bar{\pp},\tt,2}$ is a faithful exact functor.
\end{prop}
\begin{proof}
The exactness follows from \refle{lemma7.7}. Every map in $\mathcal{C}_{\bar{\pp},\tt,2}$ factors into a composition of surjection followed by an injection. \refle{lemma7.7} implies that $R^1\Gamma_{\kk,\tt}(\cdot)$ maps a nonzero surjection to a nonzero surjection and a nonzero injection to a nonzero injection.
\end{proof}
\begin{prop}
Suppose $M$ is a $(\gg,\kk)$-module. Then  $\Ext^i_{\gg,\kk}(M,R^1\Gamma_{\kk,\tt}(N))\cong\Ext^{i+1}_{\gg,\tt}(M,N)$ for $i\geq 0$.
\end{prop}
\begin{proof}
Apply the Frobenius Reciprocity Spectral Sequence in Ch. 6 of \cite{V2}, then quote \refle{lemma7.7}.
\end{proof}
\begin{corollary}
a) If $M$ is a finite dimensional $\gg$-module, then $\Ext^i_{\gg,\kk}(M,R^1\Gamma_{\kk,\tt}(N))\cong\Ext^{i+1}_{\gg,\tt}(M,N)$.

b) Suppose $N_1$ and $N_2$ are objects in $\mathcal{C}_{\bar{\pp},\tt,2}$. Then, $\Hom_{\gg,\kk}(R^1\Gamma_{\kk,\tt}(N_1), R^1\Gamma_{\kk,\tt}(N_2))\cong\Ext^1_{\gg,\tt}(R^1\Gamma_{\kk,\tt}(N_1),N_2)$ as finite-dimensional vector spaces.Thus, $\dim\Hom_\gg(N_1,N_2)\leq\dim\Ext^1_{\gg,\tt}(R^1\Gamma_{\kk,\tt}(N_1),N_2)$.
\end{corollary}

Assume again that $\tilde{\kk}$ is symmetric. Let $E$ be a simple finite dimensional $\mm$-module. Then $R^1 \Gamma_{\kk,\tt}(N_\pp(E))$ is a $(\gg,\kk)$-module of finite type, and hence a $(\gg,\tilde{\kk})$-module of finite type over $\tilde{\kk}$, i. e. a Harish-Chandra module. By the Comparison Principle [PZ4,Proposition 2.6], we have a $\gg$-module isomorphism $R^1 \Gamma_{\kk,\tt}(N_\pp(E))\cong R^1 \Gamma_{\tilde{\kk},\tt\oplus C(\kk)}(N_\pp(E))$.

The $(\gg,\kk)$-module $R^1\Gamma_{\tilde{\kk},\tt\oplus C(\kk)}(N_\pp(E))$, denoted by $A(\pp,E)$, has been studied extensively in the Harish-Chandra module literature. (See for example \cite{KV}.)
\begin{corollary}
If $\mu\geq\frac{1}{2}\lambda_1$, the Harish-Chandra module $A(\pp,E)$ has a simple socle, and $\Soc A(\pp,E)\cong R^1 \Gamma_{\kk,\tt}(L_\pp(E))$.
\end{corollary}

\section{Six examples}
In this section we consider six different pairs $(\gg,\kk)$ such that $\rk\gg=2$ and $\kk\simeq sl(2)$.

\subsection{Background on the principal series of Harish-Chandra modules} 

We start by recalling the construction of the algebraic principal series of $(\gg,\ss)$-modules for a symmetric subalgebra $\ss\subset\gg$, \cite{D}. We use this construction in subsections 8.3 - 8.6 below. Let $\ss\subset\gg$ be a symmetric subalgebra of $\gg$. Denote by $\aa_I$ a maximal toral subalgebra of $\ss^\perp$. If $\ss$ is proper, $\aa_I$ is non-zero. Let $\hh_I$ be a Cartan subalgebra of $\gg$ such that $\hh_I=(\hh_I\cap\ss)\oplus\aa_I$. Choose an element $a\in\aa_I$ such that the eigenvalues of $a$ on $\gg$ are real and $C(a)=(C(a)\cap\ss)\oplus\aa_I$. Let $\pp_{I,a}=\bigoplus_{\alpha(a)\geq 0} \gg^\alpha=\mm_I\crplus\nn_I$. 

The following results are proved in \cite{D}.
\begin{prop}
a) $\gg=\ss+\pp_{I,a}$; $\ss\cap\pp_{I,a}=\mm_I=C(\aa_I)$.

b) If $\bb_I$ is a Borel subalgebra of $\mm_I$ such that $\hh_I\subset\bb_I$, then $\bb_I\crplus\nn_I$ is a Borel subalgebra of $\gg$. Hence, $\pp_{I,a}$ is a parabolic subalgebra of $\gg$.

c) If $a'\in\aa_I$ such that $C(a')=C(a)$, then $\pp_{I,a'}$ is conjugate to $\pp_{I,a}$ under the connected algebraic subgroup $S\subset\mathrm{Aut}\gg$ whose Lie algebra is $\ss$.
\end{prop}
We define an element $a\in\ss^\perp$ to be \textit{nondegenerate} if $C(a)\cap\ss^\perp$ is a toral subalgebra of $\ss^\perp$. Moreover, an \textit{Iwasawa parabolic subalgebra for the pair} $(\gg,\ss)$ is any subalgebra of the form $\pp_{I,a}$ for some nondegenerate element $a\in\ss^\perp$, such that ad$a$ has real eigenvalues in $\gg$.

Fix an Iwasawa parabolic subalgebra $\pp_I\subset\gg$. Let $L$ be a finite-dimensional simple module over $\mm_I$. Endow $L$ with a $\pp_I$-module structure by setting $\nn_I\cdot L=0$.
\begin{defi}
a) The \textit{Iwasawa principal series module corresponding to the pair} $(\pp_I,L)$ is the $(\gg,\ss)$-module $$X(\pp_I,L)=\Gamma_\ss(\Hom_{U(\pp_I)}(U(\gg),L)).$$

b) A degenerate principal series module corresponding to the pair $(\pp_I,L)$ is the $(\gg,\ss)$-module $$Y(\qq,L)=\Gamma_\ss(\Hom_{U(\qq)}(U(\gg),L)),$$ where $\qq$ is a subalgebra containing $\pp_I$ and the $\pp_I$-module structure of $L$ extends to a $\qq$-module structure.
\end{defi}
\begin{lemma}
There is an isomorphism of $\ss$-modules $$X(\pp_I,L)\cong\Gamma_\ss(\Hom_{U(\mm_I\cap\ss)}(U(\ss),L)).$$ In particular, $X(\pp_I,L)$ is a $(\gg,\ss)$-module of finite type; if $V$ is a simple finite-dimensional $\ss$-module, then $\Hom_\ss(V,X(\pp_I,L))\cong\Hom_{\mm_I\cap\ss}(V,L)$, hence $$\dim\Hom_\ss(V,X(\pp_I,L))\leq\dim V.$$
\end{lemma}
A similar statement holds for $Y(\qq,L)$.
\begin{theo}
(Harish-Chandra's subquotient theorem)
Let $M$ be a simple $(\gg,\ss)$-module. Then there exists a simple finite-dimensional $\mm$-module $L$ such that $M$ is a subquotient of $X(\pp_I,L)$.
\end{theo}
\begin{corollary}
For any simple $(\gg,\ss)$-module $M$ and for any $\ss$-type $V$, $$\dim\Hom_\ss(V,M)\leq\dim V.$$
\end{corollary}

\subsection{$\gg=sl(2)\oplus sl(2)$, $\kk$ is a diagonal $sl(2)$-subalgebra} 

Let $\gg= sl(2)\oplus sl(2)$ and $\kk$ be the diagonal $sl(2)$-subalgebra of $\gg$. The subalgebra $\kk$ is regular in $\gg$. The pair $(\gg,\kk)$ is symmetric and its Harish-Chandra modules have been studied for over half century, see \cite{GN}, \cite{B} and \cite{HC}.

The parabolic subalgebra $\pp$ is a Borel subalgebra, and $\lambda_1=\lambda_2=2$.  We have $\rho_{\nn}=2$, hence a minimal $\kk$-type $V(\mu)$ is generic if $\mu\geq\rho_\nn -1=1$ and there is a bijection between the 1-dimensional complex family $\{\nu\in\hh^*\; |\; \nu(h)=\mu -2\}$ and the set of isomorphism classes of $(\gg,\kk)$-modules with minimal $\kk$-type $V(\mu)$. Hence any simple $(\gg,\kk)$-module $M$ with minimal $\kk$-type $\mu\geq 1$ is strongly reconstructible by Theorem 4 in [PZ2]. On the other hand, \refcor{cor6.5} above, implies this fact under the stronger assumption $\mu\geq 2$. Note that for each $\mu$, there exists a 1-dimensional complex family of simple $(\gg,\kk)$-modules with minimal $\kk$-type $V(\mu)$. These modules are multiplicity-free; a self-contained purely algebraic description of these modules in given in \cite{PS2}.

We now consider the case $\mu=0$.
\begin{prop}\label{pr8.5}
For any infinite-dimensional simple $(\gg,\kk)$-module $M$ with minimal $\kk$-type $\field{C}=V(0)$ (i. e. spherical simple $(\gg,\kk)$-module), there exists an $\hh$-module $E$ such that $$M\simeq F^1(\pp,E).$$
\end{prop}
\begin{proof}
As a $\kk$-module, $M$ is isomorphic to $\bigoplus_{j\in \field{Z}_{\geq 0}}V(2j)$, and there is no finite-dimensional simple $\gg$-module with the same central character as $M$ (see for instance \cite{PS2}). 

Now choose a 1-dimensional $\hh$-module $E$ (in the case we consider, $\mm=\hh$) such that $\omega=-2$ and $F^1(\pp,E)$ has the same central character as $M$. Then $F^0(\pp,E)=0$, since otherwise $F^0(\pp,E)$ would be a finite-dimensional $(\gg,\kk)$-module with the same central character as $M$. By an application of the Euler characteristic principle [PZ1, Theorem 11], $F^1(\pp,E)$ is isomorphic as a $\kk$-module to $\bigoplus_{j\in \field{Z}_{\geq 0}}V(2j)$. Therefore, $F^1(\pp,E)$ has the same central character and the same $\kk$-character as $M$, i.e. $M\cong F^1(\pp,E)$.
\end{proof}

Note that modules $M$ as in \refprop{pr8.5} are not strongly reconstructible. Indeed, if the central character of $M$ is regular, it is not difficult to show that there are precisely two 1-dimensional modules $E_1$ and $E_2$ such that $M\simeq F^1(\pp,E_i)$ for $i=1,2$. In addition, in this case $H^1(\nn,M)^{\omega=-2}\cong E_1\oplus E_2$. In the case of a singular central character $H^1(\nn,M)^{\omega=-2}$ is a non-trivial self-extension of the unique 1-dimensional $\mm$-module $E$ such that $M\simeq F^1(\pp,E)$.

Finally, it is true that any simple $(\gg,\kk)$-module for the pair considered is weakly reconstructible. We also remark that for any $\mu\geq 0$, $F^1(\pp,E)$ can be either irreducible or reducible.

\subsection{$\gg=sl(3)$, $\kk$ is a root $sl(2)$-subalgebra} 

Let $\gg=sl(3)$ and $\kk$ be the $sl(2)$-subalgebra of $\gg$ generated by the root spaces $\gg^{\pm(\varepsilon_1-\varepsilon_2)}$. The subalgebra $\kk$ is regular in $\gg$ and $\tilde{\kk}=\kk\oplus C(\kk)$ is a symmetric subalgebra of $\gg$ isomorphic to $gl(2)$.

The parabolic subalgebra $\pp$ is a Borel subalgebra with roots $\varepsilon_1-\varepsilon_3,\; \varepsilon_3-\varepsilon_2$ and $\varepsilon_1-\varepsilon_2$. Hence $\tilde{\rho}_\nn=\varepsilon_1-\varepsilon_2,\; \rho_\nn=2$, and any simple $(\gg,\kk)$-module $M$ of finite type over $\kk$ is strongly reconstructible for $\mu\geq\rho_\nn -1\geq 1$ by Theorem 4 in \cite{PZ2}.

On the other hand, $\lambda_1=2,\; \lambda_2=1$, hence \refcor{cor6.3} above implies the strong reconstructibility under the stronger assumption $\mu\geq\frac{3}{2}$. For completeness we note that for a $\kk$-type $V(\tilde{\mu})$ of $\tilde{\kk}$, a necessary, but not sufficient condition for $V(\tilde{\mu})$ to be generic is that $\tilde{\mu}(h)\geq 1$.

Next, $\rho_\nn^\perp=\frac{1}{2}((\varepsilon_1-\varepsilon_3)+(\varepsilon_3-\varepsilon_2))(h)=\frac{1}{2}(\varepsilon_1-\varepsilon_2)(h),$ and hence $2\rho_\nn^\perp=2$.

Fix $\mu\in\field{Z}_{\geq 1}$. As $\kk$ is regular in $\gg$, by Theorem 4 of \cite{PZ2} there exists a bijection between isomorphism classes of simple $(\gg,\kk)$-modules with lowest $\kk$-type $\mu$ and $\hh$-weights $\nu$ such that $\nu(h)=\mu-2$. If $k$ is a generator of $C(\kk)$, observe that $\nu(k)$ is a free continuous parameter of $\nu$.

An Iwasawa parabolic subalgebra $\pp_I\subset\gg$ relative to $\tilde{\kk}\subset\gg$ is a Borel subalgebra of $\gg$. Hence a finite-dimensional simple $\pp_I$-module $L$ is 1-dimensional. Write $\pp_I=\hh_I\crplus\nn_I$, and $L=L_\chi$ for $\chi\in\hh_I^*$.
\begin{prop}\label{pr8.6}
The principal series module $X(\pp_I,L_\chi)$ has finite type over $\kk$ and has lowest $\kk$-type $\field{C}=V(0)$.
\end{prop}
\begin{proof}
By Frobenius reciprocity for the principal series,
\begin{equation}
\Hom_{\tilde{\kk}}(V(\tilde{\mu}),X(\pp_I,L_\chi))\cong\Hom_{\tilde{\kk}\cap\pp_I}(V(\tilde{\mu}),L_\chi).
\end{equation}

The right hand side can be computed explicitly. Write $\tilde{\tt}=\tt+Z(\tilde{\kk})$, $\tilde{\tt}$ being a Cartan subalgebra of $\tilde{\kk}$ with basis $h$ and $k$. Let $\zeta\in\tilde{\tt}^*$ satisfy $\zeta(h)=0,\; \zeta(k)=1$. Then $\{\rho,\zeta\}$ is the basis of $\tilde{\tt}^*$ dual to the basis $\{h,k\}$ of $\tilde{\tt}$. If $V(\tilde{\mu})$ is a $\tilde{\kk}$-type we can now write $\tilde{\mu}=a\rho+b\zeta$, with $a\in\field{Z}_{\geq 0}$ and $b\in\field{C}$.

Next, $\tilde{\kk}\cap\pp_I$ is a toral subalgebra of $\tilde{\tt}$ and is spanned over $\field{C}$ by $h_I:=3h+k$. The eigenvalues of $h_I$ in $V(\tilde{\mu})$ are $3a+b-6j$ for $j\in\field{Z}_{\geq 0}$, $0\leq j\leq a$, all of multiplicity one.

The single eigenvalue of $h_I$ in $L_\chi$ is $\chi(h_I)$. Hence, $\Hom_{\field{C}h_I}(V(\tilde{\mu}),L_\chi)\neq 0$ precisely when there exists $j\in\field{N}$ with $0\leq j\leq a$ such that $3a+b-6j=\chi(h_I)\in\field{C}$. Thus, by Frobenius reciprocity, $V(\tilde{\mu})$ is a $\tilde{\kk}$-type of $X(\pp_I,L_\chi)$ iff there exists $j\in\field{Z}_{\geq 0}$ with $0\leq j\leq a$ such that $b=\chi(h_I)-3a+6j$. As a consequence, if $V(\tilde{\mu})$ is a $\tilde{\kk}$-type of $X(\pp_I,L_\chi)$, then $\chi(h_I)-3a\leq b\leq\chi(h_I)+3a$.

If we restrict the action on $X(\pp_I,L_\chi)$ from $\tilde{\kk}$ to $\kk$ we see that the multiplicity of $V(a)$ in $X(\pp_I,L_\chi)$ is $a+1=\dim V(a)$. In Figure 1 we indicate the convex hull of the $\tilde{\kk}$-support of $X(\pp_I,L_\chi)$. In general, $\chi(h_I)\in\field{C}$, but in the figure we take $\chi(h_I)>0$.

\begin{figure}[here]    
  \begin{center}
     \scalebox{0.45}{\includegraphics{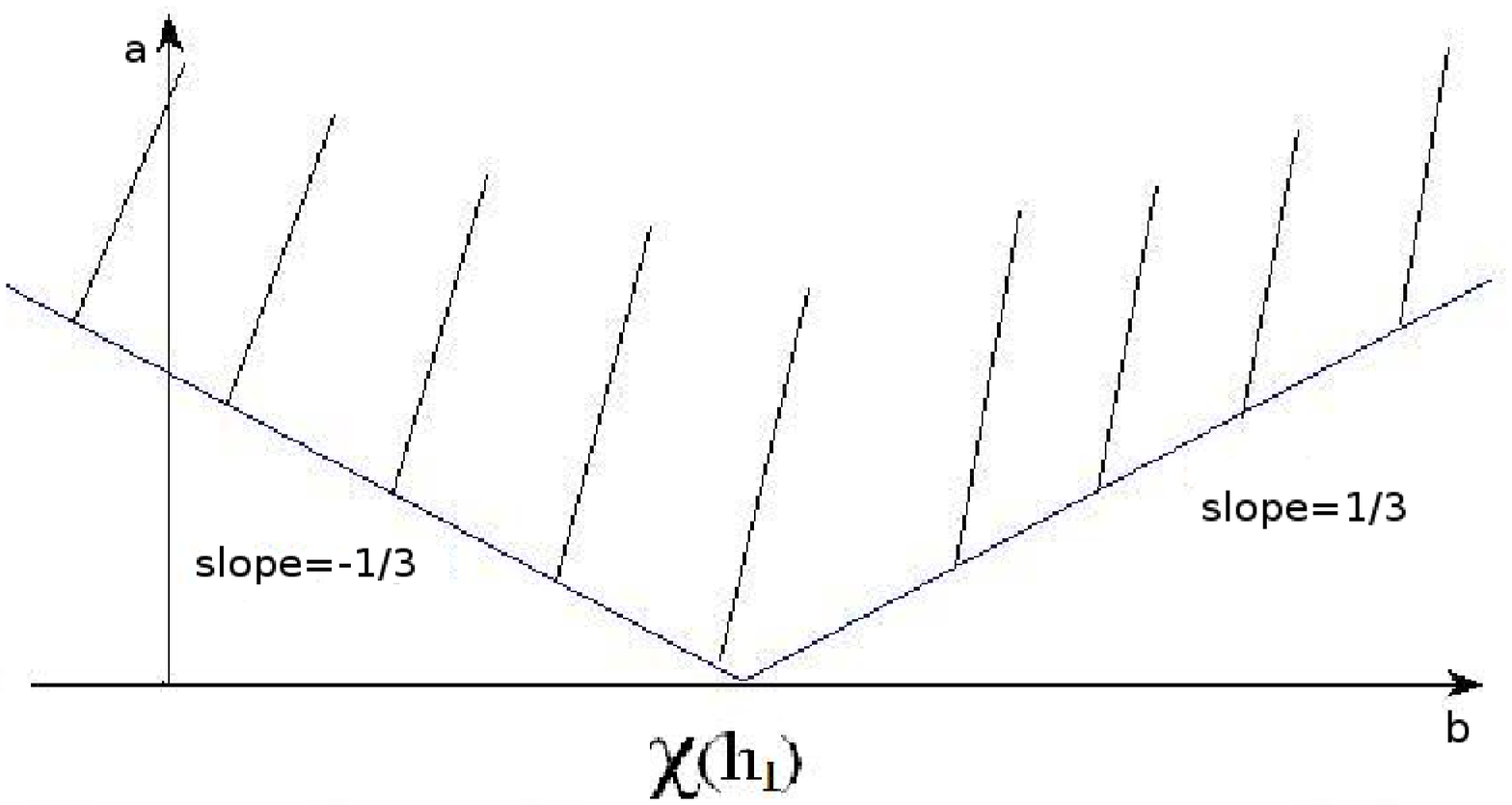}}
  \end{center}
  \caption{Figure 1}  
\end{figure}  

\end{proof}
\begin{prop}\label{pr8.7}
There exists an open dense subset $\mathcal{U}\subset\hh_I^*$ such that $X(\pp_I,L_\chi)$ is simple and not weakly reconstructible for every $\chi\in\mathcal{U}$.
\end{prop}
\begin{proof}
The existence of an open dense subset $\mathcal{U'}\subset\hh_{I}^*$ such that $X(\pp_I,L_\chi)$ is simple for $\chi\in \mathcal{U'}$ is established in \cite{Kra}. Moreover, this implies the claim as the set of weakly reconstructible modules depends on one complex and one integer parameters, while the set of irreducible principal series modules depends on two complex parameters.
\end{proof}
\begin{corollary} 
The bound $\mu\geq 1$ is sharp relative to weak (and also strong) reconstruction for $(\gg,\kk)$-modules of finite type over $\kk$.
\end{corollary}

For a classification of simple $(\gg,\tilde{\kk})$-modules, see \cite{Kra}.

\subsection{$\gg=sl(3)$, $\kk$ is a principal $sl(2)$-subalgebra} 

Let $\gg=sl(3)$ and $\kk=so(3)$, the principal $sl(2)$-subalgebra of $\gg$. The subalgebra $\kk$ is regular in $\gg$ and it is a symmetric subalgebra of $\gg$. The parabolic subalgebra $\pp$ is a Borel subalgebra and has positive roots $\varepsilon_1-\varepsilon_2$, $\varepsilon_2-\varepsilon_3$ and $\varepsilon_1-\varepsilon_3$. Hence, $\tilde{\rho}_\nn=\varepsilon_1-\varepsilon_3$. Moreover, it is easy to check that $\rho_\nn =4$, which shows that $\mu\in\field{Z}_{>0}$ is generic if $\mu\geq\rho_\nn -1=3$. On the other hand, $\lambda_1=4,\; \lambda_2=2$, so the condition $\mu\geq\frac{1}{2}(\lambda_1+\lambda_2)$ is equivalent to the same inequality: $\mu\geq\frac{1}{2}(4+2)=3$. Furthermore, we have $2\rho_\nn^\perp=4$.

Fix $\mu\in\field{Z}_{\geq 3}$. By Theorem 4 of \cite{PZ2}, there exists a bijection between isomorphism classes of simple $(\gg,\kk)$-modules with lowest $\kk$-type $\mu$ and $\hh$-weights $\nu$ such that $\nu(h)=\mu-4$. If $k$ is a generator of $\kk^\perp\cap\hh$, observe that $\nu(k)$ is a free continuous parameter of $\nu$.

Let $\pp_I\subset\gg$ be an Iwasawa parabolic subalgebra relative to $\kk$. The principal series $X(\pp_I,L_\chi)$ has two free complex parameters. Let $I_\chi$ be the sum in $X(\pp_I,L_\chi)$ of the $\kk$-types $V(0), V(2), V(4),...$. Since $\gg\cong\kk\oplus V(4)$, it is easy to see that $I_\chi$ is a $\gg$-submodule of $X(\pp_I,L_\chi)$ (see for instance [V2, Ch. 4]). A much deeper fact is that $I_\chi$ splits as a direct summand of two submodules $J_\chi$ and $K_\chi$, where the lowest $\kk$-type of $J_\chi$ is $0$ and the lowest $\kk$-type of $K_\chi$ is $2$. Furthermore, we have the following.

\begin{prop}
There exists an open dense subset $\mathcal{U}\subset\hh_I^*$ such that $K_\chi$, for $\chi\in\mathcal{U}$, is simple and not weakly reconstructible.
\end{prop}
\begin{proof}
The existence of an open dense subset $\mathcal{U'}\subset\hh_I^*$, such that the modules $J_\chi$ and $K_\chi$ are simple for $\chi\in \mathcal{U'}$, is established in [V2, Ch. 8]. This implies the claim as (similarly to the proof of \refprop{pr8.7}) the set of weakly reconstructible modules depends on one complex and one integer parameters, while the set of irreducible principal series modules depends on two complex parameters.
\end{proof}
\begin{corollary}
The bound $\mu\geq 3$ is sharp relative to weak (and also strong) reconstruction for $(\gg,\kk)$-modules of finite type over $\kk$.
\end{corollary}

\subsection{$\gg=sp(4)$, $\kk$ is a long root $sl(2)$-subalgebra} 

Let $\gg=sp(4)$. The $\hh$-roots of $\gg$ are $\pm 2\varepsilon_1, \pm 2\varepsilon_2, \pm(\varepsilon_1-\varepsilon_2),\pm(\varepsilon_1+\varepsilon_2)\in\hh^*$. Let $\kk$ be the $sl(2)$-subalgebra generated by $\gg^{\pm 2\varepsilon_1}$. The nilradical of $\pp$ has roots $\varepsilon_1-\varepsilon_2,\varepsilon_1+\varepsilon_2,2\varepsilon_1$. Hence, $\tilde{\rho}_\nn=2\varepsilon_1$, $\rho_\nn=2$, and a weight $\mu\geq 0$ is generic if $\mu\geq\rho_\nn-1=1$. On the other hand, $\lambda_1=2,\lambda_2=1$, so the condition $\mu\geq\frac{1}{2}(2+1)\geq\frac{3}{2}$ is stronger than the genericity condition. Note that $2\rho_\nn^\perp=2$. Finally, $\mm=\tt\oplus C(\kk)$, where $C(\kk)$ is the $sl(2)$-subalgebra generated by $\gg^{\pm 2\varepsilon_2}$.

Fix $\mu\in\field{Z}_{\geq 1}$. By Theorem 3 of \cite{PZ2}, we have a bijection between the following sets:

a) isomorphism classes of simple $(\gg,\kk)$-modules $M$ having finite type over $\kk$ and lowest $\kk$-type $\mu$;

b) isomorphism classes of simple, finite-dimensional $\mm$-modules $E$ such that the highest weight $\nu$ of $E$ satisfies $\nu(h)=\mu-2$.

Moreover, if $k$ is a generator of $\hh\cap C(\kk), \nu(k)$ is a free but discrete parameter of $\nu$; in this the fundamental series of $(\gg,\kk)$-modules depends on two discrete parameters.

Let $\pp_I\subset\gg$ be an Iwasawa parabolic subalgebra relative to the symmetric subalgebra $\tilde{\kk}=\kk\oplus C(\kk)$. We can choose a Levi decomposition $\pp_I=\mm_I\crplus\nn_I$ such that $\mm_I\cap\tilde{\kk}$ is the diagonal $sl(2)$-subalgebra in $\tilde{\kk}$. Let $L_\chi$ be a simple finite-dimensional $\mm_I$-module with highest weight $\chi\in\hh_I^*$, where $\hh_I$ is a Cartan subalgebra of $\mm_I$. The weight $\chi$ has one discrete and one continuous parameter.
\begin{prop}
The principal series module $X(\pp_I,L_\chi)$ has finite type over $\kk$ and has lowest $\kk$-type $0$.
\end{prop}
\begin{proof}
It is completely analogous to the proof of \refprop{pr8.6}, and Figure 2 is the analogue of Figure 1.

\begin{figure}[here]    
  \begin{center}
     \scalebox{0.45}{\includegraphics{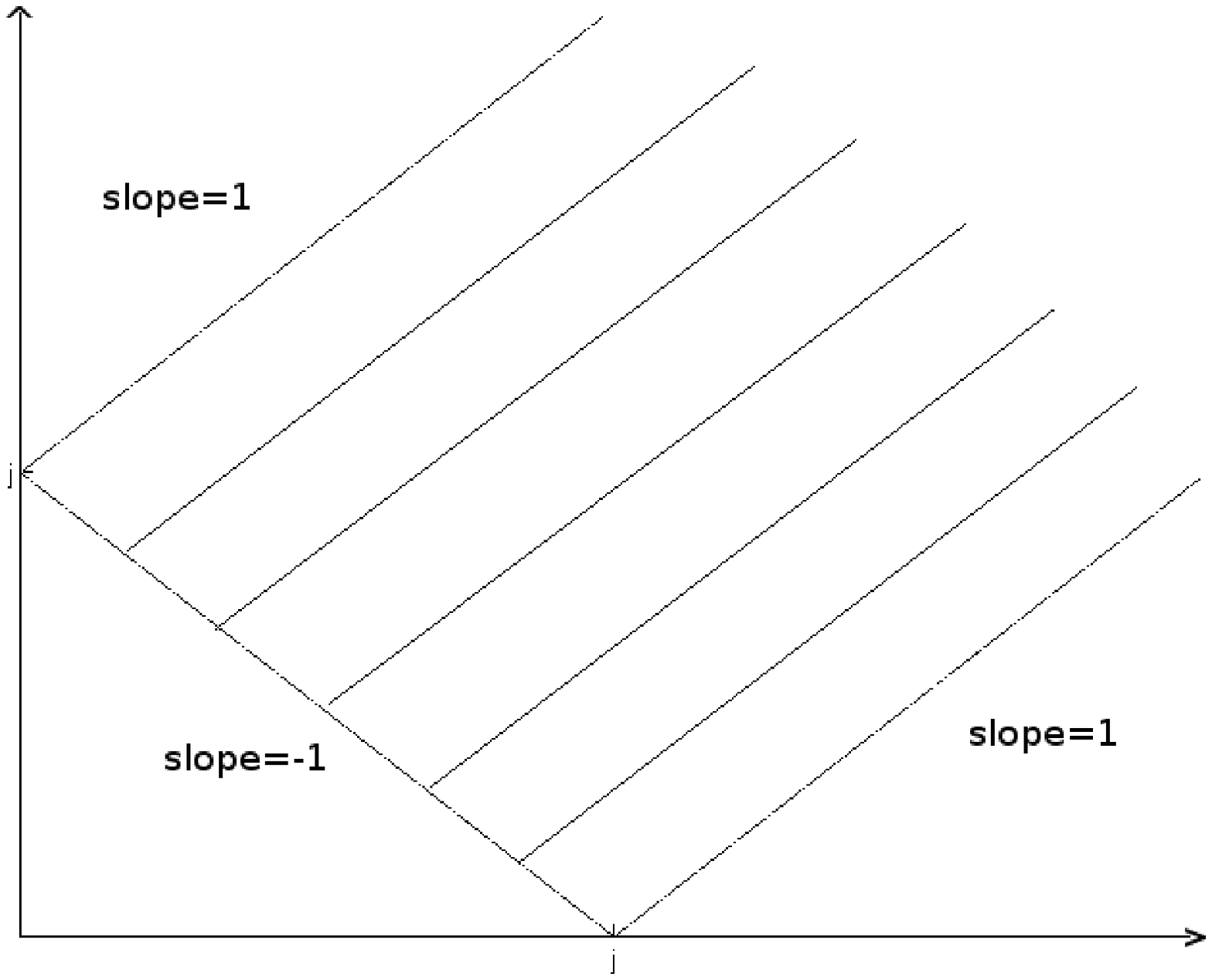}}
  \end{center}
  \caption{Figure 2}  
\end{figure} 

\end{proof}
\begin{corollary}
Every simple $(\gg,\tilde{\kk})$-module has finite type over $\kk$.
\end{corollary}
\begin{prop}
If $\chi$ is nonintegral as a weight of $\gg$, then $X(\pp_I,L_\chi)$ is simple and not weakly reconstructible.
\end{prop}
\begin{proof}
See \cite{Co}, Theorem 2.3.1.
\end{proof}
\begin{corollary} 
The bound $\mu\geq 1$ is sharp relative to weak (and also strong) reconstruction for $(\gg,\kk)$-modules of finite type over $\kk$.
\end{corollary}
\begin{proof}
The set of simple modules of the form $X(\pp_I,L_\chi)$ is not countable while the set of simple weakly reconstructible $(\gg,\kk)$-modules is countable. This implies the claim.
\end{proof}

\subsection{$\gg=sp(4)$, $\kk$ is a short root $sl(2)$-subalgebra} 

Let $\gg=sp(4)$ and $\kk$ be generated by $\gg^{\pm(\varepsilon_1-\varepsilon_2)}$. Then $\varepsilon_1(h)=1,\varepsilon_2(h)=-1$. The nilradical of the parabolic subalgebra $\pp$ has roots $\varepsilon_1-\varepsilon_2,2\varepsilon_1$ and $-2\varepsilon_2$. Hence, $\tilde{\rho}_\nn=\frac{3}{2}(\varepsilon_1-\varepsilon_2)$, $\rho_\nn=3$, and $\mu\in\field{Z}_{\geq 0}$ is generic if $\mu\geq\rho_\nn(h)-1=2$. On the other hand, $\lambda_1=2,\lambda_2=2$, so the condition $\mu\geq\frac{1}{2}(\lambda_1+\lambda_2)\geq 2$ is equivalent to being generic. Note that $2\rho_\nn^\perp=4$. Finally, $\mm=\tt\oplus C(\kk)$, where $C(\kk)$ is the $sl(2)$-subalgebra generated by $\gg^{\pm(\varepsilon_1+\varepsilon_2)}$.

Fix $\mu\in\field{Z}_{\geq 2}$. Theorem 3 of \cite{PZ2}, or equivalently \refcor{cor6.5},  implies that we have a bijection between the following sets:

a) isomorphism classes of simple $(\gg,\kk)$-modules $M$ having finite type over $\kk$ and lowest $\kk$-type $\mu$;

b) isomorphism classes of simple finite-dimensional $\mm$-modules $E$ such that the highest weight $\nu$ of $E$ satisfies $\nu(h)=\mu-4$. 

If $h'$ is a generator of $[\gg^{\varepsilon_1+\varepsilon_2},\gg^{-\varepsilon_1-\varepsilon_2}]$, observe that $\nu(h')$ is a free but discrete parameter for $\nu$.

We will exhibit a simple $(\gg,\kk)$-module $M$ of finite type over $\kk$ such that $M$ has lowest $\kk$-type $1$ but $M$ is not weakly reconstructible. Let $\tilde{\kk}=\kk\oplus C(\kk)\cong\mathrm{gl}(2)$. Let $\pp_I$ be an Iwasawa parabolic subalgebra of $\gg$ relative to $\tilde{\kk}$. This is a Borel subalgebra. Let $\qq$ be a maximal parabolic subalgebra of $\gg$ such that $\qq\supset\pp_I,\qq\neq\pp_I.$ (There are two choices for $\qq$.) Write $\qq=\ll\crplus\uu$, where $\ll$ is a reductive part of $\qq$ and $\uu$ is the nilradical of $\qq$. Note that we can choose $\ll$ so that $\ll\cap\kk=0$ and $\ll\cap\tilde{\kk}$ is a 1-dimensional toral subalgebra of $\ll$.

Next, let $L_\chi$ be a simple finite-dimensional $\ll$-module with $\hh_I$-highest weight $\chi$. Write $Y(\qq,L_\chi)$ for the degenerate principal series module $\Gamma_{\tilde{\kk}}(\Hom_{U(\qq)}(U(\gg),L_\chi))$. Since $\gg\cong\tilde{\kk}\oplus 2V(2)$ as a $\kk$-module, $Y(\qq,L_\chi)$ is a direct sum of two submodules, $Y(\qq,L_\chi)_0$ and $Y(\qq,L_\chi)_1$ corresponding to even highest weights of $\kk$ and odd highest weights of $\kk$, respectively.

\begin{lemma}
a) $Y(\qq,L_\chi)$ is a $(\gg,\kk)$-module of finite type over $\kk$.

b) The lowest $\kk$-type of $Y(\qq,L_\chi)_0$ is $\field{C}=V(0)$; the lowest $\kk$-type of $Y(\qq,L_\chi)_1$ is $V(1)$.

c) We can choose $\chi$ so that the central character of $Y(\qq,L_\chi)$ is not equal to the central character of a fundamental series module for $(\gg,\kk)$.
\end{lemma}
\begin{proof}
Straightforward calculation.
\end{proof}
\begin{prop}\label{pr8.14}
The bound $\mu\geq 2$ is sharp relative to weak reconstruction for $(\gg,\kk)$-modules of finite type over $\kk$.
\end{prop}
\begin{proof}
We take $M$ to be a simple quotient of $U(\gg)\cdot(Y(\qq,L_\chi)_0[1])$ and we chose $\chi$ so that $M$ does not have the central character of a fundamental series module.
\end{proof}

\subsection{$\gg=sp(4)$, $\kk$ is a principal $sl(2)$-subalgebra} 

Let $\gg=sp(4)$ and let $\kk$ be a principal $sl(2)$-subalgebra. Here $\mm=\hh$. It is easy to check that $\rho_\nn=7$, and the results of \cite{PZ2} (see formula (16) in \cite{PZ2}) imply that any simple $(\gg,\kk)$-module with minimal $\kk$-type $V(\mu)$ for $\mu\geq 6$ is strongly reconstructible (and in particular is of finite type). The same follows from \refcor{cor6.6} under the weaker assumption that $\mu\geq 5$.

Since $\rho_\nn^\perp=6,\lambda_1=6,\lambda_2=4$, \refcor{cor6.6} implies that, for any $\mu\in\field{Z}_{\geq 5}$, we have a bijection between the set $\{\nu\in\hh^*\; |\; \nu(h)=\mu-12\}$ and the set of isomorphism classes of $(\gg,\kk)$-modules with minimal $\kk$-type $V(\mu)$.

\begin{prop}
The bound $\mu\geq 5$ is sharp relative to the theorem on strong reconstruction for $(\gg,\kk)$-modules of finite type over $\kk$.
\end{prop}
\begin{proof}
In \cite{PS2} a simple multiplicity-free $(\gg,\kk)$-module $M_0$ with $\kk$-character $V(4)\oplus V(10)\oplus V(16)\oplus...$ and central character $\theta_{M_0}=\theta_{\frac{3}{2}e_1+\frac{1}{2}e_2}$ is exhibited: see equation 6.2 in \cite{PS2}. On the other hand, there are 8 fundamental series modules of the form $F^1(\bb,E)$ such that $\theta_{F^1(\bb,E)}=\theta_{\frac{3}{2}\varepsilon_1+\varepsilon_2}$. A non-difficult computation shows that their respective minimal $\kk$-types are $V(10),V(9),V(8),V(5),V(5),V(2),V(1)$ and $V(0)$. This shows that $M_0$ is not strongly reconstructible.
\end{proof}

We do not know whether the bound $\mu\geq 5$ is sharp relative to weak reconstruction. The following examples  demonstrate that \refth{th7.4} does not extend to the case $0\leq\mu <\frac{\lambda_1}{2}$.

\textbf{Example 1.} There is a unique 1-dimensional $\bb$-module $E_0$ such that $\theta_{F^1(\bb,E_0)}=\theta_{\frac{3}{2}\varepsilon_1+\varepsilon_2}$ and such that the minimal $\kk$-type of $F^1(\bb,E_0)$ is $V(0)$. By direct computation, $X_0=R^1\Gamma_{\kk,\tt}(L_\bb(E_0))$ is multiplicity free over $\kk$. By comparison with the simple multiplicity free modules discussed in \cite{PS2}, we conclude that $\mathrm{ch}_\kk X_0$ is the sum of two simple characters. By the Duality Theorem and the fact that $X_0$ is multiplicity free, we conclude that $X_0$ is the direct sum of two simple submodules with lowest $\kk$-types $V(0)$ and $V(4)$ respectively.

This decomposition is consistent with \refcon{conj7.4}. Moreover, the proper inclusions of $(\gg,\kk)$-modules $\bar{F}^1(\bb,E_0)\subset R^1\Gamma_{\kk,\tt}(L_\bb(E_0))\subset F^1(\bb,E_0)$ demonstrate that the inclusions discussed in \refprop{pr7.3} b) are generally proper. 

\textbf{Example 2.} There is a unique 1-dimensional $\bb$-module $E_1$ such that  $\theta_{F^1(\bb,E_1)}=\theta_{\frac{3}{2}\varepsilon_1+\varepsilon_2}$ and the minimal $\kk$-type of $F^1(\bb,E_1)$ is $V(1)$. As in Example 1, we find that $R^1\Gamma_{\kk,\tt}(L_\bb(E_1))$ is a direct sum of two simple multiplicity free $(\gg,\kk)$-modules with lowest $\kk$-types $V(1)$ and $V(3)$ respectively.

\section{Towards an equivalence of categories}
Recall the categories $\mathcal{C}_{\bar{\pp},\tt,n}$ and $\mathcal{C}_{\kk,n}$ introduced in Section 7.
\refprop{pr7.2} and \refprop{pr7.8} imply that $R^1\Gamma_{\kk,\tt}$ is a well-defined faithful and exact functor between $\mathcal{C}_{\bar{\pp},\tt,n+2}$ and $\mathcal{C}_{\kk,n}$ for $n\geq 0$.
\begin{conjecture}\label{con9.1}
Let $n\geq\frac{1}{2}(\lambda_1+\lambda_2)$. Then $R^1\Gamma_{\kk,\tt}$ is an equivalence between the categories $\mathcal{C}_{\bar{\pp},\tt,n+2}$ and $\mathcal{C}_{\kk,n}$.
\end{conjecture}

\refth{th7.4} implies that if $n\geq\frac{\lambda_1}{2}$, the simple objects $L(E)$ of the category $\mathcal{C}_{\bar{\pp},\tt,n+2}$ are being mapped by $R^1\Gamma_{\kk,\tt}$ into simple objects of $\mathcal{C}_{\kk,n}$, and \refcor{cor6.4} ensures that, under the stronger condition $n\geq\frac{\lambda_1+\lambda_2}{2}$, $R^1\Gamma_{\kk,\tt}$ induces a bijection on the isomorphism classes of simple objects of $\mathcal{C}_{\bar{\pp},\tt,n+2}$ and $\mathcal{C}_{\kk,n}$.

 \refcon{con9.1} implies the existence of an isomorphism $$\Ext^1_{\mathcal{C}_{\bar{\pp},\tt,n+2}}(L(E_1),L(E_2))\cong\Ext^1_\gg(R^1\Gamma_{\kk,\tt}(L(E_1)),R^1\Gamma_{\kk,\tt}(L(E_2)))$$ for any simple objects $L(E_1),L(E_2)$ of $\mathcal{C}_{\bar{\pp},\tt,n+2}$ where $n\geq\frac{1}{2}(\lambda_1+\lambda_2)$. We have checked the existence of such an isomorphism by direct computations in the cases of subsections 8.2 and 8.4.

In conclusion we note that it is easy to check that \refcon{con9.1} holds for the case when $\kk=\gg\cong sl(2)$. In this case $\frac{1}{2}(\lambda_1+\lambda_2)=0$ and $R^1\Gamma_{\kk,\tt}$ is an equivalence of the categories $\mathcal{C}_{\bar{\pp},\tt,n+2}$ and $\mathcal{C}_{\kk,n}$ for any $n\geq 0$.

\end{document}